%% file: main.tex
\newif\ifarXiv
\arXivtrue

\ifarXiv
	\documentclass[10pt,onecolumn,twoside]{IEEEtran}
\else
	\documentclass[10pt,twocolumn,twoside]{IEEEtran}
\fi
\input{preamble}

\acrodef{ADMM}		[ADMM]			{Alternating Direction Method of Multipliers}
\acrodef{MC}		[MC]			{Monte Carlo}
\acrodef{NR}		[NR]			{Newton-Raphson}
\acrodef{NRC}		[NRC]			{Newton-Raphson Consensus}
\acrodef{aNRC}		[a-NRC]			{asynchronous Newton-Raphson Consensus}
\acrodef{raNRC}		[ra-NRC]		{robust asynchronous Newton-Raphson Consensus}
\acrodef{raAC}		[ra-AC]			{robust asynchronous Average Consensus}

\newcommand{\flagtr}[1]{\mathtt{flag}_{\mathtt{transmission},#1}}
\newcommand{\flagrec}[1]{\mathtt{flag}_{\mathtt{reception},#1}}
\newcommand{\flagup}[1]{\mathtt{flag}_{\mathtt{update},#1}}

\title{Newton-Raphson Consensus under asynchronous and lossy communications for
  peer-to-peer networks \thanks{This work is partially supported by the Celtic
    Plus project \emph{SENDATE-Extend} (C2015/3-3), the Swedish research council
    Norrbottens Forskningsr{\aa}d project \emph{DISTRACT}, and by the European
    Research Council (ERC) under the European Union's Horizon 2020 research and
    innovation programme (grant agreement No 638992 - OPT4SMART).  
  }}

\author{Nicoletta Bof \and Ruggero Carli \and Giuseppe Notarstefano \and Luca Schenato \and Damiano Varagnolo%
  \thanks{N. Bof, R.\ Carli and L.\ Schenato are with the Department of Information Engineering, University of Padova, Via Gradenigo 6/a, 35131 Padova, Italy {\tt\small \{ bof | carlirug | schenato \}@dei.unipd.it}.}
  \thanks{Giuseppe Notarstefano is with the Department of Engineering, Universit\`a del Salento, Via per Monteroni, 73100 Lecce, Italy {\tt\small giuseppe.notarstefano@unisalento.it}.}
  \thanks{D.\ Varagnolo is with the Department of Computer Science, Electrical and Space Engineering, Lule{\aa} University of Technology, Forskargatan 1, 97187 Lule{\aa}, Sweden {\tt\small damiano.varagnolo@ltu.se}.}
  }

\begin{document}

\maketitle

\begin{abstract}
  In this work we study the problem of unconstrained convex-optimization in a
  fully distributed multi-agent setting which includes asynchronous computation
  and lossy communication. In particular, we extend a recently proposed
  algorithm named Newton-Raphson Consensus by integrating it with a
  broadcast-based average consensus algorithm which is robust to packet
  losses. We show via the separation of time scales principle that under mild
  conditions (i.e., persistency of the agents activation and bounded consecutive
  communication failures) the proposed algorithm is proved to be locally
  exponentially stable with respect to the optimal global solution. Finally, we
  complement the theoretical analysis with numerical simulations that are based
  on real datasets.
 \end{abstract}

\acresetall 

\section{Introduction}
\label{sec:introduction}

Recently, we have been witnessing a surge of interest in distributed
optimization, and in particular in distributed convex optimization. The reason
is twofold: the first is due to the advent of Big-Data analytics, whose problems
can be often cast as a large-scale convex optimization problems via Machine
Learning tools~\cite{Slavakis:14}. As so, parallelization of computation is
ought in order to obtain rapid solutions. The second reason is the advent of
Internet-of-Things and Smart Cyber-physical Systems, where a large multitude of
electronic devices are capable of sensing, communicating, and of autonomous
decision making through cooperation~\cite{Kyoung:12}. Even in this second
scenario, several estimation and control problems such as localization,
map-building, sensor calibration, power flow optimization can be cast as
large-scale convex optimization problems. The main difference between these two
scenarios is that in the former the bottleneck is mainly given by computation
time and therefore the typical architecture adopted is server-client (i.e.,
memory is centralized at a master node or redundant via synchronized
cloud architectures, and computation is parallelized among many nodes).
Our work will mainly focus on the second scenario; however, since the boundary
between the two is sometimes blurred, we will briefly overview the most relevant
literature on distributed convex optimization in general.

To cope with real-world requirements, distributed convex algorithms need to be
designed to work under asynchronous, directed, faulty and time-varying
communications. 

A popular class of algorithms that are able to cope with asynchronous updates
and lossy communication is the one of \emph{distributed subgradient
  methods}. They are simple to implement, can cope with non-differentiable
convex cost functions, and require only the computation of local
(sub)-gradients.  However, these algorithms exhibit sub-linear converge rates
even if the cost functions are smooth~\cite{Nedic:10,Nedic2010constrained}.
%
Recent works based on this approach have extended these results to directed and
possibly time-varying communication in both
discrete-time~\cite{lin2016distributed,Nedic2015distributed} and continuous-time
settings~\cite{gharesifard2014,kia2014}, however the use of a diminishing
step-size tacitly implies that the communication is synchronous (since the
step-size is designed as a function of the global time that triggers the
algorithm). Moreover, the underlying assumption for guaranteeing convergence is
that the transmitter nodes should know which packets are transmitted
successfully. This assumption corresponds to employing communication protocols
with reliable packet transmission-acknowledge mechanisms, which might be
difficult or expensive to implement over wireless media. The recent work
\cite{lee2016asynchronous} proposes an asynchronous algorithm, based on random
projections, in which the step-size (both diminishing and constant) is
uncoordinated among agents.

Another popular class of distributed optimization algorithms is the one of
\emph{dual decomposition schemes}. In this case the related literature is very
large and we refer to~\cite{yang2011distributed} for a comprehensive
tutorial. Among these algorithms, the \ac{ADMM} has attracted the attention of
the scientific community for its simple distributed implementation and good
convergence speed. This algorithm was originally proposed in mid '70s as a
general convex optimization strategy, then exploited in the context of networked
optimization~\cite{Schizas2008consensus}, and recently popularized by the
survey~\cite{boyd2011distributed}. Substantial research has been dedicated in
optimizing the free parameters of \ac{ADMM} in order to obtain fastest
convergence rates, but these are mainly restricted to synchronous
implementations over undirected communication
graphs~\cite{ghadimi2014optimalparameter,teixeira2013optimal,Nishihara:15,Iutzeler:16}. Some
recent exceptions extend dual decomposition, \cite{notarnicola2016asynchronous},
and \ac{ADMM}, \cite{Wei:13,Bianchi:16,Tsung:16}, to asynchronous scenarios with
edge-based or node-based activation schemes.

A third class of optimization algorithms, usually referred to as
\emph{Newton-based methods}, consists of strategies that exploit second-order
derivatives, i.e., the Hessians of the cost functions for computing descent
directions. For example in~\cite{wei2013newton1,wei2013newton2} the authors
apply quasi-Newton distributed descent schemes to general time-varying directed
graphs. Another approach, based on computing Newton-Raphson directions through
average consensus algorithms, has been proposed
in~\cite{Varagnolo2016newton}. Even if initially proposed for synchronous
implementations, this scheme has been later extended to cope with asynchronous
symmetric gossip communication schemes~\cite{zanella2012asynchronous}.

Finally, a different approach, based on the exchange of active constraints, has
been proposed in \cite{Notarstefano:11} for convex (and abstract) optimization
problems and extended in \cite{burger2014polyhedral} by means of cutting-plane
methods. The proposed algorithms work under asynchronous, directed and
unreliable communication.

Although there exists a large body of literature on distributed convex
optimization schemes employing synchronous and asynchronous communications, no
work has directly addressed situations where the communications are unreliable
and lossy. Unfortunately, trying to make the aforementioned algorithms cope with
packet losses using na{\"i}ve modifications (e.g., using the most recently
received message from the neighboring nodes, interpretable as using delayed
information in the algorithms) may destroy some of the hypotheses that guarantee
the convergence of the original algorithms (e.g., the doubly stochasticity or
the invariance of some quantities such as the global averages). Distributed
convex optimization in the presence of lossy communications is thus a
non-trivial task, and recently some works have specifically addressed this
problem in \ac{ADMM} schemes~\cite{Zhang:14,Peng:16,Tsung:16}; however these
strategies are restricted to networks with server-client communication
topologies.

\emph{The main contribution of this work is to propose a Newton-based algorithm
  which is robust to both asynchronous updates and packet losses and which is
  suitable for general peer-to-peer networks}. More specifically, we robustify
the Newton-Raphson approach initially proposed in~\cite{Varagnolo2016newton} by introducing a new consensus algorithm, which is an ad hoc  merging of two known schemes for consensus:
\emph{i)} the \emph{ratio} or \emph{push-sum}
consensus, useful to compute averages in networks with directed communication
graphs (i.e., networks using broadcast protocols~\cite{benezit2010weighted});
\emph{ii)} the \emph{robust consensus} algorithm, which allows for a robust
computation of arithmetic averages over networks with lossy
communication~\cite{Dominguez:11}. 
The new scheme is then able to deal with asynchronous and lossy communication protocols.
Under mild conditions, i.e., persistency of
(asynchronous) node updates, uniformly bounded consecutive communication link
failures, and connectivity of the communication graph, we then show that the
optimization algorithm is locally exponentially stable with respect to the
global solution as long as the step-size of the updates is smaller than a
certain critical value and the cost functions are sufficiently smooth. The proof
is based on time-scale separation and Lyapunov theory, and extends the results
in~\cite{ECC15}, where the convergence was proved only for quadratic cost
functions. We complement the theoretical results with numerical simulations
based on real datasets under lossy, broadcast communication. It is worth
mentioning that the algorithm we propose not only handles asynchronous updates,
as some recent references, but also is robust to packet losses.

The paper is organized as follows: \Section~\ref{sec:problem_formulation} formulates our problem and working assumptions. \Section~\ref{sec:aNRC} presents the building blocks of the scheme proposed in this manuscript. \Section~\ref{sec:raNRC} then introduces the main distributed optimization algorithm and gives some intuitions on the convergence properties of the scheme, summarized then in \Section~\ref{sec:theoretical_analysis}. Finally, \Section~\ref{sec:numerical_experiments} collects some numerical experiments corroborating the theoretical results, while \Section~\ref{sec:conclusions} draws some concluding remarks and future research directions.

\section{Problem Formulation and Assumptions}
\label{sec:problem_formulation}

We consider the separable optimization problem
\begin{equation}
	\label{eq:problem}
	x^{\ast} \DefinedAs \argmin_x f(x) = \argmin_x \sum_{i=1}^N f_i(x)
\end{equation}
where $x\in\mathbb{R}^n$ and where the local costs $f_i : \mathbb{R}^n \mapsto \mathbb{R}$ satisfy: 
\begin{assumption}[Cost smoothness]
  Each $f_i$ is known only to node $i$ and is $\C^{2}$ and strongly convex,
  i.e., its Hessian is bounded from below, $\nabla^2 f_{i}(x) > cI_n$ for all
  $x$, with $c>0$ some positive scalar\footnote{With a little abuse of notation
    we use the symbols $\nabla f(\cdot)$ and $\nabla^2 f(\cdot)$ to indicate the
    gradient and Hessian of the cost function $f(\cdot)$, respectively.}.
	\label{ass:cost_smoothness}
\end{assumption}
The communication among nodes is modeled via a communication graph that satisfies the following:
\begin{assumption}[Network connectivity]
	The communication graph among the nodes is fixed, directed and strongly connected, i.e., for each pair of nodes there is at least one directed path connecting them.
	\label{ass:network connectivity}
\end{assumption}
More formally, the communication graph is represented as $\mathcal{G} = (\V,\E)$
with nodes $\V=\left\{1,\ldots,N\right\}$ and edges $\E\subseteq \V \times \V$
so that $(i,j) \in \E$ iff node $j$ can directly receive information from node
$i$. With $\mathcal{N}_i^{\text{out}}$ we denote the set of \emph{out-neighbors}
of node $i$, i.e.,
$\mathcal{N}_i^{\text{out}} \DefinedAs \left\{j \in \V \; | \; (i,j) \in \E,
  i\neq j \right\}$
is the set of nodes receiving messages from $i$. Similarly, with
$\mathcal{N}_i^{\text{in}}$ we denote the set of \emph{in-neighbors} of $i$,
i.e.,
$\mathcal{N}_i^{\text{in}} \DefinedAs \left\{j \in \V \; | \; (j,i) \in \E\
  i\neq j \right\}$.
Their cardinality is indicated by $|\mathcal{N}_i^{\text{out}}|$ and
$|\mathcal{N}_i^{\text{in}}|$ respectively.

\begin{remark}
  In some distributed systems, as Wireless Sensor Networks, the communication
  graph is often undirected, in the sense that a node can transmit to any node
  from which it can receive. However, communication is typically only
  half-duplex, i.e., two nodes cannot communicate simultaneously, so that
  protocols with multiple communication rounds and reliable acknowledge (ACK)
  mechanisms are needed for bidirectional communication. This, in turn,
  requires pairwise synchronization and results in substantial delays; as so,
  algorithms that are suitable for broadcast-based (directed) communication
  without ACK, such as UDP, are extremely valuable also for undirected graphs.
\end{remark}

As for the concept of time, we assume that the local variables at each node are
updated at discrete time instants (e.g., based on local and possibly
unsynchronized clocks, or based on events like receiving a packet). Thus, from a
global perspective, we collect and order all time instants when at least one
variable in one node is updated and refer to it as the sequence
$\{t_k\}_{k=1}^{\infty}$. With a little abuse of notation we will then write
$x(k)=x(t_k)$ and we will study the time evolution of the nodes variables as a
discrete-time system.

Our objective is to design an algorithm solving~\eqref{eq:problem} with the following features:
\begin{enumerate}[label=F\arabic*),ref=F\arabic*]
\item \emph{Asymptotic global estimation}: each agent wants to obtain an estimate of global minimizer that asymptotically converges to the optimal solution $x^*$.
	\label{feature:asymptotic-global-estimation}
      \item \emph{Peer-to-peer (leaderless)}: each node has limited
        computational and memory resources and it is allowed to communicate
        directly only with its neighbors; moreover there is no leader/master
        node, and the communication graph is arbitrary (but strongly connected).
	\label{feature:leaderless}
      \item \emph{Distributed}: the update-rule of the local variables at each
        node depends only on the variables stored by the local node and
        by its neighbors; in other words, no multi-hop information exchange is
        allowed.
	\label{feature:distributed}
      \item \emph{Asynchronous}: events as the update of the local
        variables, and the transmission/reception of
        messages do not need to be synchronized within the node itself nor with
        its neighbors, i.e., any communication and update protocol can be used
        (e.g., time-triggered, event-triggered or hybrid). Therefore, none, one
        or multiple nodes can communicate or update their variables at any given
        time.
	\label{feature:asynchronous}
      \item \emph{Lossy broadcast communication without ACK}: communication can
        be broadcast-based with no ACK mechanisms and allow for packet losses
        (due to ambient noise, collisions, or other effects) without impairing the
        convergence properties of the algorithm.
	\label{feature:lossy-broadcast}
\end{enumerate}
To the best of authors' knowledge, none of the previously cited works possesses
all the previous features.

\section{Building blocks}
\label{sec:aNRC}

The algorithm we propose consists of three different building blocks: \emph{i)
  Newton-Raphson Consensus}, proposed in~\cite{Varagnolo2016newton} to solve
problem~\eqref{eq:problem}, \emph{ii) the push-sum algorithm}, initially
proposed in~\cite{benezit2010weighted} as an asynchronous average consensus
protocol, and \emph{iii) the robust ratio consensus algorithm}, initially
proposed in~\cite{Dominguez:11} as a robust average consensus protocol.
While possessing the first three features mentioned above
(i.e.,~\ref{feature:asymptotic-global-estimation}, \ref{feature:leaderless}, and
\ref{feature:distributed}), Newton-Raphson Consensus is nonetheless limited
since it assumes synchronous and reliable communications. The two adopted
consensus schemes (ratio consensus and its robust version) are nonetheless
limited since assume respectively reliable communications and synchronous
updates.


The major contribution of this work is to suitably modify and integrate the
three schemes above to design a distributed optimization algorithm that solves
problem~\eqref{eq:problem} and that exhibits all the
features~\ref{feature:asymptotic-global-estimation}-\ref{feature:lossy-broadcast}
above. The main challenge in doing this is that the interaction between these
algorithms might lead to instability unless some suitable assumptions are
considered. The key mathematical machinery that will be used to this means is
Lyapunov theory and separation of time-scales. 

Before providing the description of the proposed algorithm, we offer a brief
description of the three aforementioned algorithms.

\subsection{Newton-Raphson Consensus}
\label{ssec:newton-raphson-consensus}

Newton-Raphson Consensus~\cite{Varagnolo2016newton} is based on the observation that the standard Newton-Raphson update in the standard centralized scenario with a single agent can be written as
\begin{eqnarray*} 
x^+ \!\!\! &=& x - \varepsilon (\nabla^2 f(x))^{-1}\nabla f(x) \\
&=&  (1-\varepsilon)x+ \varepsilon \big(\nabla^2 f(x)\big)^{-1}\big(\nabla^2 f(x)\, x- \nabla f(x)\big)\\
&=& \!\!\! (1\!-\!\varepsilon)x\!+\! \varepsilon \big(\!\underbrace{\sum_i \!\nabla^2 f_i(x)}_{\IDefinedAs h(x)}\!\big)\!^{-1}\!\big(\!\underbrace{\sum_i (\nabla^2\! f_i(x)\, x\!-\! \nabla \!f_i(x))}_{\IDefinedAs g(x)}\big)
\end{eqnarray*}
where we used the simplified notation $x^+$ to indicate \text{$x(k+1)$} and $x$ to indicate $x(k)$. This system is exponentially stable as long as the parameter $\varepsilon>0$, which acts as a stepsize, is chosen in a proper way. If we now assume that all agents can have a different value of $x_i$ and we mimic the previous algorithm, we get the $N$ local updates:
%
%
\begin{equation} \label{eqn:NRC}
	x_i^+ \!\! =\!\!\! (1\!-\!\epsilon)x_i\!+\! \varepsilon \big(\!\underbrace{\sum_j \!\underbrace{\nabla^2 f_j(x_j)}_{\IDefinedAs h_j(x_j)}}_{\IDefinedAs \overline{h}(x_1,\ldots,x_N)}\!\big)\!^{-1}\!\big(\!\underbrace{\sum_j (\underbrace{\nabla^2\! f_j(x_j)\, x_j\!-\! \nabla \!f_j(x_j))}_{\IDefinedAs g_j(x_j)}}_{\IDefinedAs \overline{g}(x_1,\ldots,x_N)}\big) .
\end{equation}
The dynamics of the $N$ local systems is identical and exponentially stable, therefore, since they are all driven by the same forcing term $\kappa(x_1,\ldots,x_n)=(\overline h(x_1,\ldots,x_N))^{-1}\overline g(x_1,\ldots,x_N)$, intuitively we expect that
$$ x_i-x_j\to 0 , \ \ \forall i,j\;, $$
which implies that all local variable will be identical. If this is the case,
then the dynamics of each local system will eventually become the dynamics of a
standard centralized Newton-Raphson algorithm. This algorithm, however, requires
each agent to be able to instantaneously compute the two sums
$\overline h, \overline g$, which is obviously not possible in a distributed
computation set-up.  The original paper~\cite{Varagnolo2016newton} extends the
standard Newton-Raphson algorithm into a distributed scenario via the use of
synchronous lossless average consensus protocols that compute these sums
asymptotically, while~\cite{zanella2012asynchronous} extends it to the case of
asynchronous gossip-based lossless average consensus strategies.

\subsection{Push-sum Consensus}
\label{ssec:push-sum_consensus}

The Newton-Raphson Consensus scheme described
in \Section~\ref{ssec:newton-raphson-consensus} requires each node to compute
the two sums $y_i = \overline g$ and $z_i =\overline h $ at least asymptotically
in order to apply a Newton-Raphson descent. In fact, since the ratio of the two
quantities is needed, each agent can asymptotically converge to a scaled version
of the two. That is, assuming each variable $x_i$, $i \in \mathcal{V}$, to be
fixed, we require
\begin{eqnarray*} 
	y_i &\to& \eta_i \overline{g}(x_1,\ldots,x_N) = \eta_i \sum_j g_j(x_j)\\
	z_i &\to& \eta_i \overline{h}(x_1,\ldots,x_N) = \eta_i \sum_j h_j(x_j),
\end{eqnarray*}
where $\eta_1, \ldots, \eta_N$ are possibly time-dependent, non-zero
scalars. Here the right arrow means that the difference between left and right
hand-sides goes to zero as the iteration counter goes to infinity. Having identified our aim, we first describe the push-sum algorithm, which is able to solve the given problem in an asynchronous communication scenario. Then, we describe the robust ratio consensus which is able to solve the problem in a scenario where the communication is unreliable but the protocol is synchronous. One of the aim of this work will be the merging of these two schemes to obtain a robust and asynchronous consensus algorithm.

%
Under synchronous communication, the local updates of the \emph{push-sum} or \emph{ratio consensus} introduced
in~\cite{benezit2010weighted} are, for each
$i \in \mathcal{V}$,
\begin{eqnarray} 
	y^+_i &=& \frac{1}{|\mathcal{N}^{\text{out}}_i|+1}y_i + \sum_{j\in \mathcal{N}^{\text{in}}_i} \frac{1}{|\mathcal{N}^{\text{out}}_j|+1}y_j
	\label{equ:ratio-consensus-synchronous-y} \\
	z^+_i &=& \frac{1}{|\mathcal{N}^{\text{out}}_i|+1}z_i + \sum_{j\in \mathcal{N}^{\text{in}}_i} \frac{1}{|\mathcal{N}^{\text{out}}_j|+1}z_j,
	\label{equ:ratio-consensus-synchronous-z}
\end{eqnarray}
paired with the initialization $y_i(0) =
g_i(x_i)$, $z_i(0) = h_i(x_i)$.
Assuming for notation simplicity a scalar optimization problem, the previous
update can be written as
\begin{eqnarray*} 
	\bm{y}^+&=& P \bm{y}  \\
	\bm{z}^+&=& P \bm{z},
\end{eqnarray*}
where $\bm{y}=[y_1\,\cdots\,y_N]^T,\, \bm{z}=[y_1\,\cdots\,z_N]^T$. In this way the matrix $P$ results to be column-stochastic and its induced graph $\mathcal{G}_P$ (i.e., $(i,j) \in \mathcal{G}_P$ if $[P]_{ji} \neq 0$) coincides with the original communication graph (i.e., $\mathcal{G}_P = \mathcal{G}$). Since we assume $\mathcal{G}$ to be strongly connected, this guarantees that\footnote{These well-known results can also be readily derived from standard theories on Markov Chains~\cite{Seneta:06}.}
\begin{equation*}
	\begin{array}{lll}
		y_i &\to& \eta_i \sum_i y_i(0) = \eta_i \sum_i g_i(x_i) = \eta_i \overline{g}(x_1,\ldots,x_N)\\
		z_i &\to& \eta_i \sum_i z_i(0) = \eta_i \sum_i h_i(x_i) = \eta_i \overline{h}(x_1,\ldots,x_N)
	\end{array}
\end{equation*}
where $\eta = [\eta_1\,\cdots\,\eta_N]^T$ is the right eigenvector of $P$
relative to the unique unitary eigenvalue, i.e., $P\eta=\eta$ and
$\eta_i>0, \forall i$.

The ratio consensus described above can then be extended to asynchronous
implementations (as proposed in~\cite{benezit2010weighted}). Let at any time $k$
only one node $i$ activate, update its variables, and broadcast them to its
out-neighbors, and then, consistently, let the generic receiving node $j$ update
its local variables. The update rules for $y_i$ and $y_j$ therefore become
\begin{eqnarray} 
	y^+_i &=& \frac{1}{|\mathcal{N}^{\text{out}}_i|+1}y_i \label{eqn:PS_1} \\
	y^+_j &=& y_j + \frac{1}{|\mathcal{N}^{\text{out}}_i|+1}y_i = y_j +y_i^+
	\quad \forall j \in  \mathcal{N}^{\text{out}}_i
	\label{eqn:PS_2}
\end{eqnarray}
(the rules for $z_i$ and $z_{j}$ being equal in structure). In this scenario, the global dynamics can be described by a time-varying consensus matrix that depends on the specific node that is activated, i.e. $P(k)\in\{P_1,\ldots,P_N\}$, where the matrices $P_i$ are still column-stochastic. As shown via weak ergodic theory considerations in \cite{benezit2010weighted}, if the activation of the nodes is randomized and i.i.d.\ then the local variables converge to
\begin{equation}
	\begin{array}{l}
		y_i \to \eta_i(k) \sum_i y_i(0)=\eta_i(k) \sum_i g_i(x_i)=\eta_i(k) \overline{g}(x_1,\ldots,x_N)\\
		z_i \to \eta_i(k) \sum_i z_i(0)=\eta_i(k) \sum_i h_i(x_i)=\eta_i(k) \overline{h}(x_1,\ldots,x_N)
	\end{array}
\label{eqn:ratio_limit}
\end{equation}
where $\eta_i(k)>0$ is time-varying and depends on the activation sequence of the nodes.

\subsection{Robust ratio consensus}
\label{ssec:robust_ratio_consensus}

The synchronous ratio-consensus strategy defined by iterations~\eqref{equ:ratio-consensus-synchronous-y} and \eqref{equ:ratio-consensus-synchronous-z} in \Section~\eqref{ssec:push-sum_consensus} loses its convergence properties in case of lossy communications. A na{\"i}ve attempt to solve this problem is then to use a buffer such that when $i$ does not receive a message from $j$ then $i$ updates its local variables by using the latest values that it has received from $j$. Focusing only on~\eqref{equ:ratio-consensus-synchronous-y} to avoid repetitions, mathematically this corresponds to add an additional local variable $y_{i}^{(j)}$ representing the latest $y_{j}$ received by $i$ from $j$, and to transform the update rule~\eqref{equ:ratio-consensus-synchronous-y} into
\begin{eqnarray*} 
	y^{(j)+}_i &=&
	\left\{
		\begin{array}{ll}
			y_j & \text{if $y_j$ is received} \\
			y^{(i)}_j & \text{otherwise}
		\end{array} \ \ \    \forall i \in \mathcal{V} , \forall j \in \mathcal{N}^{\text{in}}_i
	\right. \\
	y^+_i &=&
	\frac{1}{|\mathcal{N}^{\text{out}}_i|+1} y_i
	+
	\sum_{j \in \mathcal{N}^{\text{in}}_i}
	\frac{1}{|\mathcal{N}^{\text{out}}_j|+1}y^{(j)}_i ,\ \ \forall i \in \mathcal{V}. \\
\end{eqnarray*}
However this solution does not preserve the total mass of the variables $y_i$
during the progress of the algorithm, i.e.,
\text{$\sum_i y_i(k) \neq \sum_i y_i(0)$}, differently from the original
lossless ratio consensus; this eventually leads the average consensus algorithm
not to converge to the desired value.
 
To overcome this issue, it is possible to add some additional ``mass counter''
variables $\sigma_{i,y},\rho^{(j)}_{i,y}$ that guarantee the preservation of the
masses even in the presence of packet losses~\cite{Dominguez:11}. More
specifically, in this way the synchronous
update~\eqref{equ:ratio-consensus-synchronous-y} transforms into
\begin{eqnarray} 
\sigma^+_{i,y} &=& \sigma_{i,y}+y_i, \ \  \forall i \in \mathcal{V}  \label{eqn:RRC_1}\\
\rho^{(j)+}_{i,y} &=& \left\{ \begin{array}{ll} \sigma_{j,y} &  \mbox{If $\sigma_{j,y}$ is received} \\
\rho^{(j)}_{i,y} & \mbox{otherwise}  \end{array} \ \ \    \forall j \in \mathcal{N}^{\text{in}}_i \right. \label{eqn:RRC_2}\\
y^+_{i} &=&\!\! \frac{1}{|\mathcal{N}^{\text{out}}_i|+1}y_{i}\! +\! \! \sum_{j\in \mathcal{N}^{\text{in}}_i} \frac{1}{|\mathcal{N}^{\text{out}}_j|+1}(\rho^{(j)+}_{i,y}\!\! -\! \rho^{(j)}_{i,y}) \ \ \ \ \label{eqn:RRC_3}
\end{eqnarray}%
where the ``mass counter'' variables are initialized to zero, i.e., $\sigma_{i,y}(0) = \rho^{(j)}_{i,y}(0)=0$ for every $i$ and $j$. 

Observe that, each node $i$ has a counter $\sigma_{i,y}(k)$ to keep track of the
total $y$-mass sent to its neighbors up to iteration $k$, and, for each
neighbor $j\in\mathcal{N}_i^{\text{in}}$, a counter $\rho^{(j)}_{i,y}(k)$ to
take into account the total $y$-mass received from $j$ up to iteration $k$. If
during iteration $k$ node $i$ receives information from node $j$, the
information related to node $j$ used in the update of the variable $y_i$ is
$\nu^{(j)}_{i,y}(k) = \sigma_{j,y}-\rho^{(j)}_{i,y}$.
%
%
The fictitious variable $\nu^{(j)}_{i,y}(k)$ corresponds to a ``virtual mass"
stored on edge $(j,i) \in \mathcal{E}$. Under reliable transmission, such
virtual mass is zero, while each time a packet loss occurs, this variable
accumulates the additional mass that node $j$ wants to transfer to node $i$, and
therefore it is not lost. As so, the total mass stored on the nodes and the
edges is preserved, regardless of the packet loss sequence.
Thus, for each time instant $k$,
\begin{equation}
	\sum_i
	\left(
		y_i(k)
		+
		\sum_{j\in \mathcal{N}^{\text{in}}_i}
	  \nu^{(j)}_{i,y}(k) 
	\right)
	=
	\sum_i y_i(0).
	\label{eqn:invariant}
\end{equation}

Let $\bm{y}$ and $\bm{\nu}_y$ be the vectors collecting, respectively, the
variables $y_i$, $i \in V$, and $\nu^{(j)}_{i,y}$, $i \in V$ and
$j \in \mathcal{N}^{\text{out}}_i$, and, accordingly, let $\bm{y}_a$ be the
augmented variable defined as
$\bm{y}_a = \left[\bm{y}^T \,\,\bm{\nu}_y^T\right]^T$. Similarly let
$\bm{z}_a = \left[\bm{z}^T \,\,\bm{\nu}_z^T\right]^T$, it can be shown that
$$
\bm{y}_a^+ =  M(k) \bm{y}_a, \qquad \bm{z}_a^+ =  M(k) \bm{z}_a
$$
where $M(k)$ is an augmented column-stochastic matrix, and, from weak ergodicity
theory, that local variables $y_i$, $z_i$, converge asymptotically as
in~\eqref{eqn:ratio_limit}~\cite{Dominguez:11}. {As it will be clear in the
  next sections, matrix $M(k)$ will be a building block for the design and
  analysis of our distributed optimization algorithm.}

\section{The \acf{raNRC}}
\label{sec:raNRC}

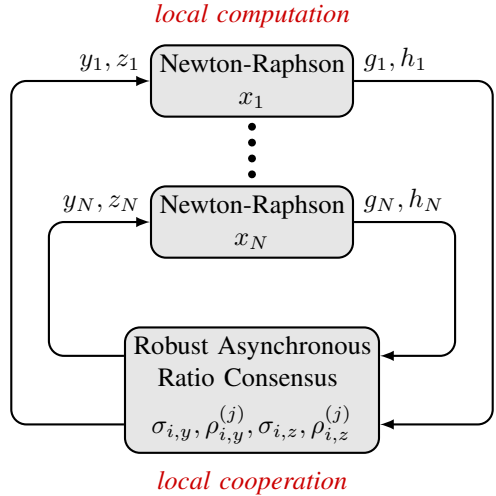
\begin{figure}[!thbp]
	\centering
	\input{figure_NR_lossy_D}
\caption{Graphical representation of the \acf{raNRC}.}
\label{fig:Figure_NR_lossy}
\end{figure}

This section merges the three building blocks \emph{Newton-Raphson Consensus},
\emph{push-sum consensus} and \emph{robust ratio consensus} into one algorithm,
called \acf{raNRC}, that solves problem~\eqref{eq:problem} and exhibits all the
features listed in \Section~\ref{sec:problem_formulation}. The algorithm can be
organized in a block scheme as in \Figure~\ref{fig:Figure_NR_lossy}.

We propose a ``meta distributed algorithm'' which can result in different
distributed algorithms depending on the (possibly asynchronous and packet-lossy)
communication protocol implemented in the network.
The meta algorithm consists of four main blocks of code implemented by each node
$i\in\mathcal{V}$ in the network: \emph{Initialization} (at startup),
\emph{Data Transmission}, \emph{Data Reception} and \emph{Estimate Update}.

Except for the first block, which corresponds to a one-time execution at
startup, the blocks can be executed asynchronously, with possibly different
execution rates. The scheduling of these three blocks, for each agent $i$, is determined by three
binary variables $\flagtr{i},\flagrec{i},\flagup{i}$
whose evolutions are determined by the communication protocol.
Each code block is assumed to be executed \emph{sequentially} and \emph{atomically}, i.e., the
local variables and flags cannot be changed by any other process. For example, if node is executing Estimate Update and a new packet is incoming, this packet is either dropped or placed in a buffer till Estimate Update is not completed. 
Thus, a distributed algorithm will be simply the combination of the given meta
scheme with a communication protocol defining how the flags are activated. For
example, in an event-triggered communication protocol the reception of a packet
may sequentially trigger (if no other block is being executed) the Data
Reception block, which then triggers the Estimate Update block, and that finally
triggers the Data Transmission block. In the following we assume that when an agent is idle, it is always ready to receive a new packet and when a packet is received by the $i$-th node then $\flagrec{i}$ is set to one.


One of the strengths of the proposed algorithm, is that it is independent of the
specific communication protocol as long as it satisfies some mild assumptions in
terms of minimum scheduling rate of each block and maximum consecutive
packet losses, which will be formally stated in the next section.
We are, then, ready to provide a pseudo-code description of ra-NRC as in
Algorithm~\ref{alg:raNRC}.
Notice that the local variables in the algorithm mimic the variable names and
purpose of the ones defined in the previous section.


\begin{algorithm}
\caption{\acf{raNRC} for node $i$}
    \begin{algorithmic}[1]
   	 \REQUIRE  $x^o$, $\varepsilon$, $c$
	 \Statex{ \bf Initialization} (atomic)
		\STATE $x_i \leftarrow x^o$
		\STATE $y_i \leftarrow 0$, \ $g_i \leftarrow 0$, \ $g^{old}_i \leftarrow 0$
		\STATE $z_i \leftarrow I_n$, \ $h_i \leftarrow I_n$, \ $h^{old}_i \leftarrow I_n$
		\STATE $\sigma_{i,y} \leftarrow 0$, \ $\sigma_{i,z} \leftarrow 0$
		\STATE $\rho_{i,y}^{(j)} \leftarrow 0$,  \ $\rho_{i,z}^{(j)} \leftarrow 0, \ \ \ \forall j \in \mathcal{N}^{\text{in}}_i$
		\STATE $\flagrec{i}\leftarrow 0$, \ $\flagup{i}\leftarrow 0$
		\STATE $\flagtr{i}\leftarrow 1$

\Statex{ \bf Data Transmission}	(atomic)	
		\IF {$\flagtr{i}=1$}
		\STATE $\mathtt{transmitter\_node\_ID} \leftarrow i$
		\STATE $y_i \leftarrow \frac{1}{|\mathcal{N}^{\text{out}}_i|+1}y_i$
		\STATE $z_i \leftarrow \frac{1}{|\mathcal{N}^{\text{out}}_i|+1}z_i$
		\STATE $\sigma_{i,y}  \leftarrow \sigma_{i,y}+y_i$
		\STATE $\sigma_{i,z}  \leftarrow \sigma_{i,z}+z_i$
		\STATE Broadcast: $\mathtt{transmitter\_node\_ID},\sigma_{i,y},\sigma_{i,z}$
		\STATE $\flagtr{i}\leftarrow 0$
		\ENDIF

\Statex{ \bf Data Reception}	(atomic)	
		\IF {$\flagrec{i}=1$}
		\STATE $j\leftarrow \mathtt{transmitter\_node\_ID}, \ \ \ (j \in \mathcal{N}^{\text{in}}_i)$   
		\STATE $y_i  \leftarrow y_i + \sigma_{j,y}-\rho_{i,y}^{(j)}$
		\STATE $z_i  \leftarrow z_i+ \sigma_{j,z}-\rho_{i,z}^{(j)}$
		\STATE $\rho_{i,y}^{(j)} \leftarrow \sigma_{j,y}, \ \ \ \forall j \in \mathcal{N}_i$
		\STATE $\rho_{i,z}^{(j)} \leftarrow \sigma_{j,z}, \ \ \ \forall j \in \mathcal{N}_i$
		\STATE $\flagrec{i} \leftarrow 0$
		\STATE $\flagup{i} \leftarrow 1$ (optional)
		\ENDIF

\Statex{ \bf Estimate Update}		(atomic)
		\IF {$\flagup{i}=1$}
		\STATE $x_i  \leftarrow (1-\varepsilon) x_i + \varepsilon {\color{blue} {z_i}^{-1}}  y_i$
		\STATE $g_i^{\text{old}} \leftarrow g_i$
		\STATE $h_i^{\text{old}} \leftarrow h_i$
		\STATE $h_i \leftarrow \nabla^2 f_i(x_i)$
		\STATE $g_i \leftarrow h_i x_i-\nabla f_i(x_i)$	
		\STATE $y_i  \leftarrow  y_i + g_i-g_i^{\text{old}}$
		\STATE $z_i  \leftarrow  z_i + h_i-h_i^{\text{old}}$
		\STATE $\flagup{i}\leftarrow 0$
		\STATE $\flagtr{i}\leftarrow 1$ (optional)
		\ENDIF
\end{algorithmic}
\label{alg:raNRC}
\end{algorithm}

We now provide a detailed explanation of the pseudo-code.

The first block \emph{Initialization}~(lines 1-7) is a one-time operation
preformed by each node at the beginning of the algorithm. The only free
parameter to set is the initial estimate $x^o$ for the global optimization,
while all other variables are set to zero or to identity matrices of the proper dimension.

The blocks \emph{Data Transmission}~(lines 8-15) and \emph{Data
  Reception}~(lines 17-24) implement a new Robust Asynchronous Ratio
Consensus (see bottom block in Figure~\ref{fig:Figure_NR_lossy}), which merges
the benefits of the push-sum algorithm, with its asynchronous nature, and the
robust ratio consensus with its resilience to packet losses. Moreover, our
proposed Robust Asynchronous Ratio Consensus has the advantage to be \emph{fully
  parallel}, in the sense that multiple nodes can transmit at the same time,
since any potential collision will result in a packet loss already handled by
the algorithm.  Specifically, the update of variables $y_i, z_i$ at the
time of transmission (line 10-11) are the same as in the push-sum consensus
given by Eqn.~\eqref{eqn:PS_1}. The update for $\sigma_{i,y}$ in the algorithm
(line 12) is identical to Eqn.~\eqref{eqn:RRC_1}, however the
variable $\sigma_{i,y}$ in our algorithm is based on the value of $y_i$ that has
been updated above (line 10). Therefore, variables $\sigma_{i,y}$'s in
Algorithm~\ref{alg:raNRC} are scaled by a factor
$\frac{1}{|\mathcal{N}^{\text{out}}_i|+1}$ as compared to those in
Eqn.~\eqref{eqn:RRC_1}. Since the variables $\rho_{j,y}^{(i)}$ will be just
(possibly delayed) copies of the variable $\sigma_{i,y}$ (line 21), also these
variables are scaled by a factor $\frac{1}{|\mathcal{N}^{\text{out}}_i|+1}$ as
compared to those appearing in Eqn.~\eqref{eqn:RRC_2}. Similar arguments apply
for the variables related to $\sigma_{i,z},\rho_{j,z}^{(i)}$. Once the update of
the variables has been completed, the transmitting node broadcasts only the
variables $\sigma_{i,y},\sigma_{i,z}$ and its ID to its neighbors. After
transmission, the node returns to an idle-mode (line 15). If a neighboring node
$i$ is in the receiving mode and actually receives a message (line 17), then it extracts the transmitter node ID $j$
and the corresponding variables $\sigma_{j,y},\sigma_{i,z}$ (line 18). The
variable $y_i$ is updated similarly to Eqn.~\eqref{eqn:PS_2}, where $y^+_i$ is
replaced by the term $\sigma_{j,y}-\rho_{i,y}^{(j)}$, which is the same as the
last term appearing in Eqn.~\eqref{eqn:RRC_3}\footnote{Note that since the
  packet is received, $\rho_{i,y}^{(j)+}=\sigma_{j,y}$.}. The local variable
$\rho_{i,y}^{(j)}$ is then updated (line 21) similarly to Eqn.~\eqref{eqn:RRC_2}
in the robust ratio consensus.

The last block \emph{Estimate Update} is responsible for implementing a local
version of Newton-Raphson. The update of the local estimate $x_i$ of the global
optimizer, available at each node $i$, is performed via the Newton-Raphson
Consensus described in the previous section. In practice, the roles of $y_i$ and
$z_i$ are those of (scaled) local approximations of the global functions
$\overline{g}(x_1,\ldots,x_N)$ and $\overline{h}(x_1,\ldots,x_N)$ defined
above. As so, mimicking Eqn.~\eqref{eqn:NRC}, the proposed algorithm uses these
variables to implement an approximated Newton-Raphson (line 27), where the
operator $[\cdot]_c$ is used to avoid divide-by-zero numerical issues during the
transient, and the variable $\varepsilon$ corresponds to the stepsize. Since the
local variables $x_i$ are continuously updated, also the global functions
$\overline{g}(x_1,\ldots,x_N)=\sum_i g_i(x_i)$ and
$\overline{h}(x_1,\ldots,x_N)=\sum_i h_i(x_i)$ need to be updated
accordingly. This cannot be done instantaneously due to the networked nature of
the framework and has be achieved through the asynchronous robust ratio consensus
(see Figure~\ref{fig:Figure_NR_lossy}). In order to be able to track the
continuously changing signals $g_i$ and $h_i$, each node has to compute these
signals before and after updating the $x_i$ ($g_i^{\text{old}}$ e
$h_i^{\text{old}}$ in lines (28-29) and $g_i$ e $h_i$ lines (30-31),
respectively) and then update the ``consensus" variables $y_i$ and $z_i$ in
order to track the current sums $\overline{g}(x_1,\ldots,x_N)$ and
$\overline{h}(x_1,\ldots,x_N)$ (lines 32-33). In fact, this operation guarantees
that, similarly to Eqn.~\eqref{eqn:invariant}, the following invariant are
preserved:
\begin{eqnarray}
 \sum_i ( y_i(k) + \sum_{j\in \mathcal{N}^{\text{in}}_i} (\sigma_{j,y}(k) - \rho^{(j)}_{i,y}(k) ))\!\!&=&\!\!\sum_i g_i(k), \ \ \  \label{eqn:invariant_ARNRC}\\
 \sum_i ( z_i(k) + \sum_{j\in \mathcal{N}^{\text{in}}_i} (\sigma_{j,z}(k) - \rho^{(j)}_{i,z}(k) ))\!\!&=&\!\!\sum_i h_i(k),
\end{eqnarray}
where, with a slight abuse of notation, with $g_i(k)$ and $h_i(k)$ we denote $g_i(x_i(k))$ and $h_i(x_i(k))$ respectively.
The intuition behind the convergence of the algorithm, is that if the local
estimates $x_i$ change slower than the rate at which the asynchronous robust
ratio consensus converges, which can be achieved by choosing a sufficiently
small stepsize $\varepsilon$, then we would expect that
\begin{eqnarray}
	y_i(k) &\to & \eta_i(k) \sum_i g_i(k), \ \
                      \  \label{eqn:invariant_ARNRC2}\\ 
	z_i(k) &\to & \eta_i(k) \sum_i h_i(k)             
\end{eqnarray}

A formal proof of the ra-NRC algorithm and the necessary conditions in terms of
node activation and packet loss frequencies, when a particular communication
protocol is adopted, are given in the next section.

\begin{remark}
  The \acl{raNRC} has the demanding requirements that full matrices
  $\sigma_{i,z}$ needs to be transmitted and inverted, which could be rather
  demanding if the feature space dimension $n$ is large. Similarly to what has
  been proposed in \cite{Varagnolo2016newton}, it is possible to modify the
  proposed algorithm to use Jacobi or Gradient descents which have reduced
  communication and computational requirements. More specifically, the only
  modification needed is to substitute line~(37) with the following ones
\begin{eqnarray*}
 h_i &\leftarrow& \mathrm{diag}(\nabla^2 f_i(x_i)), \ \ \text{Jacobi Descent Consensus} \\
 h_i &\leftarrow& I_n,  \hspace{2.1cm} \text{Gradient Descent Consensus}, 
 \end{eqnarray*}  
 where the operator $ \mathrm{diag}(A)$ returns a diagonal matrix whose diagonal
 elements coincide with the diagonal elements of $A$. As so, for the Jacobi
 Descent Consensus it is necessary to invert $n$ scalars and to transmit only
 the $n$ diagonal elements, while for the Gradient consensus only one scalar
 needs to be transmitted and inverted.  Of course, the price to pay with these
 choices is a likely slower convergence rate.
\end{remark}

\GN{
\begin{remark}
The results presented in this work are concerned with the local stability and guarantee that $z_i(k)>cI_n$ for all times $k$, however we have observed that in order to increase the basin of attraction and the robustness of the algorithm it is suitable to force $z_i(k)\geq cI_n$. A simple solution is to replace line~(27) with 
$$x_i  \leftarrow (1-\varepsilon) x_i + \varepsilon {\cmax{z_i}^{-1}}  y_i$$
where the operator $\cmax{\cdot}$ is defined as
\begin{equation*}
	\cmax{ z }
	\DefinedAs
	\left\lbrace
	\begin{array}{ll}
		z & \text{ if } z \geq cI_n \\
		cI_n & \text{ otherwise.} 
	\end{array}
	\right.
\end{equation*}
where $z\in\mathbb{R}^{n\times n}$ is a positive semidefinite matrix, and $I_n$
is the identity matrix of dimension $n$. This does not impair the local stability analysis provided below since close to the equilibrium point we have $\cmax{z_i(k)}^{-1} = {z_i(k)}^{-1}$.
\end{remark}
}

\section{Dynamical system interpretation of ra-NRC}

In this Section, we introduce an asynchronous and lossy communication protocol
that defines the evolution of the flags in Algorithm~\ref{alg:raNRC}, in order to carry out the convergence analysis of the algorithm itself, and we will also show that this choice of communication protocol is not restricting. The protocol selected allows us to rewrite the resulting ra-NRC as a dynamical
system of the form:
%
\begin{equation*}
\left\{
\begin{array}{l}
\bm{x}(k+1)=\bm{x}(k)+\epsilon \phi (k, \bm{x}(k), \bm{\xi}(k))\\
\bm{\xi}(k+1)=\varphi (k, \bm{x}(k)), \bm{\xi}(k)),
\end{array}
\right.
\end{equation*}
where proper definitions of variables $\bm{x}$, $\bm{\xi}$ and maps $\phi$ and
$\varphi$ can be found in  Corollary \ref{cor:dynamicalSystemRewriting}. 

In particular, we focus our analysis on an \emph{asymmetric broadcast}
communication protocol subject to packet losses, which represents a widely used
communication protocol in wireless sensor networks applications. Specifically,
let $t_0, t_1, t_2, \ldots$ be an ordered sequence of time instants, i.e.,
$t_0 < t_1 < t_2< \ldots$. We assume that at each time instant one node, say
$i$, is activated. Then, node $i$ performs in order the operations in the
\emph{Estimate Update} block and in the \emph{Data Transmission} block,
broadcasting to all its out-neighbors in $\mathcal{G}$ the updated variables
$\sigma_{i,y}, \sigma_{i,z}$. The transmitted packet might be received or not by
$j \in \mathcal{N}^{\text{out}}_i$, depending whether $(i,j)$ is reliable or
not at the time of transmission. If $(i,j)$ is reliable, then node $j$ performs, in order, the operations in the
\emph{Data Reception} block, and in the \emph{Estimate Update} block. Since there
is no risk of confusion, in the following we denote $t_k$ only by the index $k$,
referring to it as the $k$-th iteration of the ra-NRC algorithm.

An algorithmic description of the \emph{asymmetric broadcast} communication
protocol with packet losses (for the ra-NRC Algorithm~\ref{alg:raNRC}) is provided in
\Algorithm~\ref{alg:raNRC-broadcast}. Here, with a slight abuse of notation, we
denote within the parentheses after the flag variables the owner of the
corresponding variable.  Moreover, in the following, without loss of generality,
we assume that, the node performing the transmission step during the $k$-th
iteration is node $i$.

\begin{algorithm}
  \caption{Asymmetric broadcast {for} ra-NRC algorithm}
    \begin{algorithmic}[1]
\vspace{0.1cm}
	 \Statex{ \bf Node $i$ is activated} 
		\STATE $\flagup{i} \leftarrow1$ \,\,(line $26$)\,\, : \,\,\hfill {\it Estimate Update}
		\STATE $\flagtr{i}\leftarrow 1$ \,(line $8$) \, : \hfill {\it Data transmission}	
\vspace{0.2cm}
\Statex{ \bf For $j \in \mathcal{N}^{\text{out}}_i$, if $(i,j)$ is reliable}		
		\STATE $\flagrec{j}\leftarrow 1$ \,\,(line $17$)\,\, : \hfill \,\, {\it Data reception}
                 \STATE $\flagup{j} \leftarrow1$ \,\,(line $26$)\,\, : \,\, \hfill {\it Estimate Update}
\end{algorithmic}
\label{alg:raNRC-broadcast}
\end{algorithm}

In order to keep the notation lighter, from now on, we restrict to the scalar
case, i.e., $x_i \in \R$ for all $i$. {Consistently we will denote, e.g.,
  $f'_i$ and $f''_i$ respectively the first and second derivatives of the
  function $f_i$.}

Next, for the sake of analysis, we provide a sequential description of the
ra-NRC algorithm when the communication protocol in Algorithm
\ref{alg:raNRC-broadcast} is adopted. Observe that, once activated, node $i$
updates $x_i$, $g_i^{\text{old}}$, $h_i^{\text{old}}$, $g_i$, $h_i$ according to
lines $27$, $28$, $29$, $30$, $31$, i.e.,
\begin{align}\label{eq:update_x_gold_g}
x_i(k+1)&=(1-\varepsilon) x_i (k)+ \varepsilon \cmax{z_i(k)}^{-1}  y_i(k)\nonumber\\
g_i^{\text{old}}(k+1)&= g_i(k) \nonumber \\
h_i^{\text{old}}(k+1)&= h_i(k) \nonumber \\
g_i(k+1) &= f''_i(x_i(k+1))x_i(k+1)-f'_i(x_i(k+1))\nonumber \\
h_i(k+1) &= f''_i(x_i(k+1)).\nonumber
\end{align}
Based on $g_i(k+1)$ and $h_i(k+1)$, the variables $y_i$ and $z_i$ are updated performing in order the steps in lines $32$, $10$, and $33$, $11$, respectively, which result in
\begin{align*}
y_i (k+1) &= \frac{1}{|\mathcal{N}^{\text{out}}_i|+1}\left(y_i(k)+  g_i(k+1) - g_i^{\text{old}}(k+1)\right)\\
z_i (k+1) &= \frac{1}{|\mathcal{N}^{\text{out}}_i|+1}\left(z_i(k)+  h_i(k+1) - h_i^{\text{old}}(k+1)\right),
\end{align*}
and, in turn, from lines $12$, $13$, we have that
\begin{align*}
\sigma_{i,y}(k+1) &=\sigma_{i,y}(k)+y_i(k+1)\\
\sigma_{i,z}(k+1) &=\sigma_{i,z}(k)+z_i(k+1).
\end{align*}
The quantities $\sigma_{i,y}(k+1)$, $\sigma_{i,z}(k+1)$ are transmitted by node
$i$ to its out-neighbors; if $(i,j)$ is reliable, then node $j$, based on the
\emph{Data Reception} packet, updates the local variables $y_j$, $z_j$,
$\rho_{j,y}^{(i)}$, $\rho_{j,z}^{(i)}$ as\footnote{As far as the variables $y_j$ and $z_j$ are concerned, to denote their updates in the \emph{Data Reception packet} we introduce the auxiliary variables $y'_j$, $z'_j$, since the overall updates of the current values of $y_j$ and $z_j$ are performed in the subsequent \emph{Data Update} packet.}
\begin{align*}
y'_j & = y_j(k)+  \sigma_{i,y}(k+1)-\rho_{j,y}^{(i)}(k)\\ 
z'_j & =  z_j(k)+\sigma_{i,z}(k+1) -\rho_{j,z}^{(i)}(k)\\
\rho_{j,y}^{(i)}(k+1) & = \sigma_{i,y}(k+1)\\
\rho_{j,z}^{(i)}(k+1)& =\sigma_{i,z}(k+1)
\end{align*}
and, subsequently, based on the \emph{Data Update} packet, updates the local variables $x_j$, $g_j^{\text{old}}$, $h_j^{\text{old}}$, $g_j$, $h_j$, $y_j$, $z_j$ as
\begin{align*}
x_j (k+1)& = (1-\varepsilon) x_j(k) + \varepsilon \frac{y_j(k)}{\cmax{z_j(k)}}\\
g_j^{\text{old}}(k+1) & =  g_j(k)\\
h_j^{\text{old}}(k+1) &= h_j(k)\\
g_j(k+1)& = f''_j(x_j(k+1))x_j(k+1)-f'_j(x_j(k+1))\\
h_j(k+1)&= f''_j(x_j(k+1))\\
y_j(k+1) & =  y'_j + g_j(k+1)-g_j^{\text{old}}(k+1)\\
z_j(k+1)  &=  z'_j + h_j(k+1)-h_j^{\text{old}}(k+1).
\end{align*}
Next, we provide a suitable vector-form description of the Asymmetric broadcast ra-NRC algorithm.

To do so, similarly to the \emph{robust ratio consensus algorithm} revisited in Section~\ref{ssec:robust_ratio_consensus}, we first need to build an \emph{augmented} network that contains all the nodes
in $V$ and also some additional virtual nodes; precisely, a virtual node for each link in $\mathcal{E}$. Let us denote the augmented network by
$\mathcal{G}_a=\left(V_a, \mathcal{E}_a\right)$, where $V_a=V \cup \mathcal{E}$
and
$$
\mathcal{E}_a=\mathcal{E} \cup \left\{\left((i,j),j\right)\,|\,(i,j) \in \mathcal{E}\right\}\cup \left\{\left(i,(i,j)\right)\,|\,(i,j) \in \mathcal{E}\right\}.
$$
Similarly to what done in Section~\ref{ssec:robust_ratio_consensus}, for each
$(i,j) \in \mathcal{E}$ we introduce the auxiliary variables
$\nu_{j,y}^{(i)}(k)$, $\nu_{j,z}^{(i)}(k)$, defined as
\begin{align*}
\nu_{j,y}^{(i)}(k) & =\sigma_{i,y}(k)-\rho_{j,y}^{(i)}(k)\\
\nu_{j,z}^{(i)}(k) &=\sigma_{i,z}(k)-\rho_{j,z}^{(i)}(k).
\end{align*}
Recall that the role of the above variables is to keep track of the transmitted
mass, which has not been received due to packet losses.  Accordingly, let
$\bm{\nu}_y$ and $\bm{\nu}_z$ be the vectors that collect, respectively, all the
variables $\nu_{j,y}^{(i)}$ and $\nu_{j,z}^{(i)}$, $i \in V$ and
$j \in \mathcal{N}^{\text{out}}_i$. Assuming that $|\mathcal{E}| = N_\mathcal{E}$, then
$\bm{\nu}_y$, $\bm{\nu}_z\, \in \,\R^{N_\mathcal{E}}$. 
Now let 
\begin{equation*}
\bm{y}  =\left[
\begin{array}{c}
y_1\\
\vdots\\
y_N
\end{array}
\right], \qquad
\bm{z}= \left[
\begin{array}{c}
z_1\\
\vdots\\
z_N
\end{array}
\right],
\end{equation*}
and, based on these vectors,
let us build the augmented vectors $\bm{y}_a$, $\bm{z}_a \, \in \,\R^{N +N_\mathcal{E}}$ as
$$
\bm{y}_a = \left[
\begin{array}{c}
\bm{y}\\
\bm{\nu}_y
\end{array}
\right], \qquad \bm{z}_a = \left[
\begin{array}{c}
\bm{z}\\
\bm{\nu}_z
\end{array}
\right].
$$
Moreover let
\begin{align*}
\bm{g} & =\left[g_1,\ldots,g_N\right]^T \\
\bm{g}^{\text{old}} & =\left[g_1^{\text{old}},\ldots,g_N^{\text{old}}\right]^T\\
\bm{f}''(\bm{x}) \bm{x} &= \left[f''_1(x_1)\,x_1, \ldots, f''_N(x_N)\,x_N\right]^T\\ 
\bm{f}'(\bm{x}) & =\left[f'_1(x_1), \ldots, f'_N(x_N)\right]^T\\
\bm{y}/{\bm{z}}&= \left[y_1/x{z_1}, \ldots,y_N/{z_N} \right]^T.
\end{align*}
Since we are considering a lossy scenario, it might happen that the packet transmitted by node $i$ is either received or not received by node $j \in \mathcal{N}^{\text{out}}_i$. For this reason, it is convenient to introduce the sets
$$
\tilde{\mathcal{N}}_i (k)= \left\{j \in \mathcal{N}^{\text{out}}_i \,\,\, \text{such that $(i,j)$ is reliable at time }k \right\},
$$ 
and, its complement on $\mathcal{N}^{\text{out}}_i$,
$$
\bar{\mathcal{N}}_i (k)= \mathcal{N}^{\text{out}}_i  \setminus \tilde{\mathcal{N}}_i (k) .
$$
To state Proposition \ref{prop:vector-form}, where we provide a vector form
description of Algorithm \ref{alg:raNRC-broadcast}, it is convenient to resort
to the following notational convention. When referring to an $N$-dimensional
vector, we assume its components to be indexed according to the nodes in $V$, while
when referring to an $N_\mathcal{E}$-dimensional vector, we assume its components to be indexed
according to the edges in $\mathcal{E}$. In particular, $e_i \in \R^N$ and
$e_{(i,j)}\in \R^{N_\mathcal{E}}$ denote the vectors with all the components equal to zero,
except, respectively, the one related to node $i$ and the one related to edge
$(i,j)$, which are equal to one; that is $e_i$, $i \in V$ and $(e_{i,j})$,
$(i,j) \in \mathcal{E}$, are the vectors of the canonical basis of,
respectively, $\R^N$ and $R^{N_\mathcal{E}}$.

\begin{proposition}\label{prop:vector-form}
  The ra-NRC algorithm with asymmetric broadcast (Algorithm~\ref{alg:raNRC} and
  Algorithm~\ref{alg:raNRC-broadcast}), can be written in vector form
  as\footnote{We observe that the matrices $S$, $S_a$ and $M$ depend on which
    node is activated, and on which edges between this node and its
    out-neighbors are reliable. In order to keep the notation lighter, we do not
    make this dependency explicit (for instance using some superscript or
    subscript); instead, we limit ourselves to emphasize only the time-varying
    nature of these matrices, i.e., just writing $S(k)$, $S_a(k)$ and $M(k)$.}
\begin{align}\label{eq:algorithm}
\bm{x}(k+1)&=\bm{x} (k)+ \varepsilon \,S(k)\, \left(\bm{p}(k)- \bm{x}(k)\right)\\
\bm{g}^{old}(k+1)&= S(k)\bm{g}(k)+(I-S(k)) \bm{g}^{\text{old}}(k) \nonumber \\
\bm{g}(k+1) &=  \bm{f}''(\bm{x}(k+1)) \bm{x}(k+1)- \bm{f}'(\bm{x}(k+1)) \nonumber \\
\bm{h}^{old}(k+1)&= S(k)\bm{h}(k)+(I-S(k)) \bm{h}^{\text{old}}(k) \nonumber \\
\bm{h}(k+1) &=  \bm{f}''(\bm{x}(k+1)) \nonumber \\
\bm{y}_a(k+1)&=M(k)\bm{y}_a(k) + \nonumber \\
& \qquad\qquad  T(k)\, \left(\bm{g}(k+1)-\bm{g}^{\text{old}}(k+1)  \right) \nonumber\\
\bm{z}_a(k+1)&=M(k) \bm{z}_a(k) +  \nonumber\\
&\qquad\qquad  T(k) \, \left(\bm{h}(k+1)-\bm{h}^{\text{old}}(k+1)  \right) \nonumber\\ 
\bm{p}(k+1)&=\frac{\bm{y}(k+1)}{{\bm{z}(k+1)}}\nonumber
\end{align}
where 
\begin{align*}
S(k) = e_ie_i^T +\!\!\!\sum_{j \in \tilde{\mathcal{N}}_i(k)} e_je_j^T\;\;\; \text{and} \;\;\;
T(k)  =
\left[
\begin{array}{c}
T_V(k) \\
T_{\mathcal{E}}(k)
\end{array}
\right],
\end{align*}
with
\begin{align*}
T_V(k) & = \frac{1}{|\mathcal{N}^{\text{out}}_i|+1} \left(  e_ie_i^T +\sum_{j \in \tilde{\mathcal{N}}_i(k)} e_je_i^T \right) + \sum_{j \in \tilde{\mathcal{N}}_i(k)} e_je_j^T\\
T_{\mathcal{E}}(k) &= \frac{1}{|\mathcal{N}^{\text{out}}_i|+1} \sum_{j \in \bar{\mathcal{N}}_i} e_{(i,j)}e_i^T
\end{align*}
and where $M(k)$ is a column stochastic matrix such that
$$
M(k)=
\left[
\begin{array}{cc}
M_{VV}(k) & M_{V \mathcal{E}}(k) \\
M_{\mathcal{E} V}(k) & M_{\mathcal{E} \mathcal{E}}(k)
\end{array}
\right]
$$
with
\begin{align*}
M_{VV}(k) &= \frac{1}{|\mathcal{N}^{\text{out}}_i|+1} \left(  e_ie_i^T +\sum_{j \in \tilde{\mathcal{N}}_i} e_je_i^T \right) + \sum_{h \neq i} e_h e_h^T\\
M_{V \mathcal{E}}(k) & = \sum_{j \in \tilde{\mathcal{N}}_i} e_je_{(i,j)}^T \\
M_{\mathcal{E} V}(k) & =\frac{1}{|\mathcal{N}^{\text{out}}_i|+1}  \sum_{j \in \bar{\mathcal{N}}_i} e_{(i,j)} e_i^T \,\,\, \left(\,=T_{\mathcal{E}}(k) \,\right)\\
M_{\mathcal{E} \mathcal{E}}(k) &= \sum_{j \in \bar{\mathcal{N}}_i} e_{(i,j)} e_{(i,j)}^T + \sum_{(r,s)\,:\, r \neq i} e_{(r,s)} e_{(r,s)}^T.
\end{align*}
\end{proposition}
\begin{IEEEproof}
We start by observing that, only nodes in $\tilde{N}_i(k) \cup \left\{i\right\}$ update the variables $x$, $g^{\text{old}}$, $g$, $h^{\text{old}}$, $h$. Moreover, observe that the matrix $S(k)$ can be seen as a \emph{selection} matrix which selects the nodes in $\tilde{N}_i(k) \cup \left\{i\right\}$. This explains the vector form of the first five equations in \eqref{eq:algorithm}.

Now, to each edge $(i,j)$, $j \in \mathcal{N}_i^{\text{out}}$, we associate the indicator function variable $X_{i,j}(k)$ as follows:
$$
X_{i,j}(k)=
\left\{
\begin{array}{l}
1, \,\,\,\text{if} \,\,\, (i,j) \,\,\text{reliable at time } k\\
0, \,\,\,\text{if} \,\,\, (i,j) \,\,\text{not reliable  at time }k.\\
\end{array}
\right.
$$
In the following, for the sake of simplicity, we consider only the update of $\bm{y}_a$ (the update of $\bm{z}_a$ is similar). 
Recall that
\begin{align}\label{eq:y_i_Update}
y_i (k+1) &= \frac{1}{|\mathcal{N}^{\text{out}}_i|+1}\left(y_i(k)+  g_i(k+1) - g_i^{\text{old}}(k+1)\right).
\end{align}


Observe that, for $j \in \mathcal{N}_i^{\text{out}}$, by using the indicator function defined above, we can write that
$$
\rho_{j,y}^{(i)}(k+1)=X_{i,j}(k) \sigma_{i,y}(k+1) + \left(1-X_{i,j}(k)\right) \rho_{j,y}^{(i)}(k).
$$
Since
$$
\nu_{j,y}^{(i)}(k)=\sigma_{i,y}(k)-\rho_{j,y}^{(i)}(k),
$$
and
\begin{align*}
\sigma_{i,y}(k+1) &=\sigma_{i,y}(k)+y_i(k+1),
\end{align*}
it follows that
\begin{align}\label{eq:nu_j_y}
\nu_{j,y}^{(i)}(k+1)&=\left(1-X_{i,j}(k)\right) \left( \nu_{j,y}^{(i)}(k)+ \right.\\
&\,\,\,\,\left. \frac{1}{|\mathcal{N}^{\text{out}}_i|+1}\left[y_i(k) + g_i(k+1)-g_i^{\text{old}}(k+1)\right] \right) \nonumber
\end{align}
and that
\begin{align*}
&y_j(k+1)=  y_j(k)+\\
& +X_{i,j}(k) \left[ y_i(k+1) + g_j(k+1)-g_j^{\text{old}}(k+1)+\nu_{j,y}^{(i)}(k)\right] \nonumber
\end{align*}
and, in turn, that
\begin{align}\label{eq:y_j_k+1}
&y_j(k+1)= y_j(k)+ X_{i,j}(k)\frac{1}{|\mathcal{N}^{\text{out}}_i|+1} y_i(k) +  \\
&\qquad\qquad+X_{i,j}(k) \left[g_j(k+1)-g_j^{\text{old}}(k+1)+\nu_{j,y}^{(i)}(k)+ \right. \nonumber \\
&\qquad\qquad \qquad\,\,\,\,\,\,\,\,\,\left. \frac{1}{|\mathcal{N}^{\text{out}}_i|+1}\left[g_i(k+1)-g_i^{\text{old}}(k+1)\right] \right].\nonumber
\end{align}
From \eqref{eq:y_i_Update} and \eqref{eq:y_j_k+1} we can write that
\begin{align*}
&\bm{y}(k+1)  = \\
& \left[\frac{1}{|\mathcal{N}^{\text{out}}_i|+1} \left(  e_i+\sum_{j \in \tilde{\mathcal{N}}_i(k)}  e_j \right) e_i^T  + \sum_{h \neq i} e_h e_h^T\right] \bm{y}(k) + \\
& +\sum_{j \in \tilde{\mathcal{N}}_i} e_je_{(i,j)}^T \,\bm{\nu} + \left[\frac{1}{|\mathcal{N}^{\text{out}}_i|+1} \left(  e_i +\sum_{j \in \tilde{\mathcal{N}}_i(k)} e_j \right)e_i^T + \right. \\
&\qquad \qquad\qquad \,\,\,\, \left.+\sum_{j \in \tilde{\mathcal{N}}_i(k)} e_je_j^T\right] \, \left(\bm{g}(k+1)-\bm{g}^{\text{old}}(k+1)  \right)\\
& = M_{VV}(k) \bm{y} + M_{V\mathcal{E}}(k) \,\bm{\nu}_y (k)+\\
&\qquad \qquad \qquad  \qquad \qquad +  T_{V}(k) \left(\bm{g}(k+1)-\bm{g}^{\text{old}}(k+1)  \right).
\end{align*}

From \eqref{eq:nu_j_y}, we have that 

\begin{align*}
&\bm{\nu}_y (k+1)=  \left[\sum_{j \in \bar{\mathcal{N}}_i} e_{(i,j)} e_{(i,j)}^T + \sum_{(r,s)\,:\, r \neq i} e_{(r,s)} e_{(r,s)}^T\right] \bm{\nu}_y+\\
&+\left[\frac{1}{|\mathcal{N}^{\text{out}}_i|+1}  \sum_{j \in \bar{\mathcal{N}}_i} e_{(i,j)} e_i^T\right] \left(\bm{y} + \bm{g}(k+1)-\bm{g}^{\text{old}}(k+1)  \right)\\
& = M_{\mathcal{E} V}(k) \bm{y}(k) + M_{\mathcal{E} \mathcal{E}}(k) \bm{\nu}_y(k) +\\
&\qquad \qquad \qquad + T_{\mathcal{E}}(k)\left(\bm{g}(k+1)-\bm{g}^{\text{old}}(k+1)  \right).
\end{align*}

The above computations explain the vector-form illustrated in equations \eqref{eq:algorithm}.

The fact that $M(k)$ is a column-stochastic matrix can be shown by simply
verifying that the sum of the elements of each column is equal to one.
\end{IEEEproof}
Observe that variables $\bm{y}_a, \bm{z}_a$ {are trajectories of a linear,
  time-varying algorithm with column-stochastic state-matrix}, driven by the
differences $\bm{g} -\bm{g}^{\text{old}}$, $\bm{h} -\bm{h}^{\text{old}}$.
  
From the previous proposition, the next fact follows directly. 
\begin{corollary}\label{cor:dynamicalSystemRewriting}
Let
$\bm{\xi}=\left[\bm{g}^T, {\bm{g}^{\text{old}}}^T,\bm{h}^T, {\bm{h}^{\text{old}}}^T, \bm{y}_a^T, \bm{z}_a^T,\bm{p}^T\right]^T$, 
then, system in \eqref{eq:algorithm} can be written
as: 
\begin{equation}\label{eq:NonLinear_System}
\left\{
\begin{array}{l}
\bm{x}(k+1)=\bm{x}(k)+\epsilon \;\phi (k, \bm{x}(k), \bm{\xi}(k))\\
\bm{\xi}(k+1)=\varphi (k, \bm{x}(k)), \bm{\xi}(k)),
\end{array}
\right.
\end{equation}
where $\bm{x} \in \R^{N}$, $\bm{\xi} \in \R^{7N+2|\mathcal{E}|}$,
$\phi \,:\, \N \times \R^{N} \times \R^{7N+2|\mathcal{E}|}\,\to\, \R^{N}$,
$\varphi \,:\, \N \times \R^{N} \times \R^{7N+2|\mathcal{E}|}\,\to\,
\R^{7N+2|\mathcal{E}|}$,
$\epsilon>0$ and where equations in \eqref{eq:algorithm} properly define the
maps $\phi$ and $\varphi$.
\end{corollary}

Finally, we characterize a \emph{mass conservation} property of system in
\eqref{eq:algorithm}, which will be useful in next section.  
\begin{lemma}\label{eq:ConservationProperty}
Consider system in \eqref{eq:algorithm}. Then, for all $k \in \N$, the following equalities hold true
\begin{align*}
\sum_{\ell=1}^N \left( y_\ell(k)\,+\, \sum_{j \in \mathcal{N}_\ell^{\text{out}}} \nu_{j,y}^{(\ell)}(k) \right) & =  \sum_{\ell=1}^N  g_\ell(k)\\
\sum_{\ell=1}^N \left( z_\ell(k)\,+\, \sum_{j \in \mathcal{N}_\ell^{\text{out}}} \nu_{j,z}^{(\ell)}(k) \right) & =  \sum_{\ell=1}^N  h_\ell(k).
\end{align*}
\end{lemma}
\begin{IEEEproof}
We provide the proof of only the first equality; the second one can be proved analogously. We proceed by induction. The property is trivially true for $k=0$. Indeed according to the Initialization block we have that $y_\ell(0)=g_\ell(0)=g_\ell^{\text{old}}(0)= \nu_{j,y}^{(\ell)}(0)=0$ for all $\ell$ and $j \in \mathcal{N}_\ell^{\text{out}}$; the fact that $g_\ell(0)=g_\ell^{\text{old}}(0)=0$ implies that also $g_\ell^{\text{old}}(1)=0$ for all $\ell$. Now, we assume the property to be true for $k$ and we show that it holds also for $k+1$. 
Without loss of generality, assume that node $i$ is activated at iteration $k$. Then we have
\begin{align*}
&\sum_{\ell=1}^N \left(y_\ell(k+1)\,+\, \sum_{j \in \mathcal{N}_\ell^{\text{out}}} \nu_{j,y}^{(\ell)}(k+1)\right)=\1^T \,\bm{y}_a(k+1)\\
& = \1^T M(k) \bm{y}_a(k) + \1^T\,T(k)\,\left(\bm{g}(k+1)-\bm{g}^{\text{old}}(k+1)  \right) \\
& = \sum_{\ell=1}^N \left( y_\ell(k)\,+\, \sum_{j \in \mathcal{N}_\ell^{\text{out}}} \nu_{j,y}^{(\ell)}(k) \right)+  g_i(k+1)-g_i^{\text{old}}(k+1) \\
& \qquad \qquad \qquad
+\sum_{j \in \mathcal{N}_i^{\text{out}}} X_{i,j}(k) \left(g_j(k+1)-g_j^{\text{old}}(k+1)\right)\\
&= \sum_{\ell=1}^N  g_\ell(k)
+\sum_{j \,\in \,\tilde{\mathcal{N}}_i(k) \,\cup \, \left\{ i \right\} } \left(g_j(k+1)-g_j^{\text{old}}(k+1)\right),
\end{align*}
where, in the above equalities, we have used the properties
$$
\1^T M(k) = \1^T, \qquad \1^T \,T(k) = e_i^T + \sum_{j \in \tilde{\mathcal{N}}_i}  e_j^T
$$
and the inductive hypothesis
$$
 \sum_{\ell=1}^N \left( y_\ell(k)\,+\, \sum_{j \in \mathcal{N}_\ell^{\text{out}}} \nu_{j,y}^{(\ell)}(k) \right) = \sum_{\ell=1}^N  g_\ell(k).
$$
By simple algebraic manipulations, we can write 
\begin{align*}
 &\sum_{\ell=1}^N  g_\ell(k) + \sum_{j \in \tilde{\mathcal{N}}_i(k) \,\cup \, \left\{ i \right\}}  \left(g_j(k+1)-g_j^{\text{old}}(k+1)\right) \\
 &= \sum_{j \in \tilde{\mathcal{N}}_i(k) \,\cup \, \left\{ i \right\}} g_j(k) + \sum_{j \notin \tilde{\mathcal{N}}_i(k) \,\cup \, \left\{ i \right\}} g_j(k) + \\
 &\qquad\qquad +\sum_{j \in \tilde{\mathcal{N}}_i(k) \,\cup \, \left\{ i \right\}}  \left(g_j(k+1)-g_j^{\text{old}}(k+1)\right)\\
 &=\sum_{j \in \tilde{\mathcal{N}}_i(k) \,\cup \, \left\{ i \right\}}  g_j(k+1) +\sum_{j \notin \tilde{\mathcal{N}}_i(k) \,\cup \, \left\{ i \right\}} g_j(k)+\\
 &\,\,\,\,\qquad\qquad+ \sum_{j \in \tilde{\mathcal{N}}_i(k) \,\cup \, \left\{ i \right\}}  \left(g_j(k)-g_j^{\text{old}}(k+1)\right).
\end{align*}
Now, observe that, if $\ell \notin \tilde{\mathcal{N}}_i(k) \,\cup \, \left\{ i \right\}$ then $g_\ell(k+1)=g_\ell(k)$, and, if $\ell \in \tilde{\mathcal{N}}_i(k) \,\cup \, \left\{ i \right\}$ then $g_\ell^{\text{old}}(k+1)=g_\ell(k)$. Then, from the previous expression, it follows
\begin{align*}
 &\sum_{\ell=1}^N  g_\ell(k) + \sum_{j \in \tilde{\mathcal{N}}_i(k) \,\cup \, \left\{ i \right\}}  \left(g_j(k+1)-g_j^{\text{old}}(k+1)\right) \\
 &= \sum_{\ell=1}^N g_\ell(k+1).
\end{align*}
This concludes the proof.
\end{IEEEproof}

\begin{remark}
In this Section, we have provided a dynamical system description of ra-NRC algorithm, assuming the asymmetric broadcast communication protocol has been adopted. However, it is worth stressing that similar computations hold also for other communication protocols like \emph{symmetric gossip}, \emph{asymmetric gossip}, \emph{coordinated broadcast}\footnote{For a concise but effective description of the aforementioned protocols we refer the interested reader to \cite{IFAC17}.}. When adopting one of the above aforementioned communication protocols it turns out that ra-NRC algorithm can again be described as in \eqref{eq:algorithm}, with the only difference related to the matrix $M(k)$ which is still a column stochastic matrix but with a slight different structure, { and to the selection matrix $S(k)$ }. This justifies the fact that the convergence results we provide in the next Section, which are specifically tailored to the scenario considered in this Section, can be technically extended to also other types of communication protocols.
\end{remark}

\section{Theoretical analysis of the ra-NRC}
\label{sec:theoretical_analysis}

We now provide a theoretical analysis of the Asymmetric broadcast ra-NRC algorithm, described in Algorithm \ref{alg:raNRC-broadcast}. In particular, we provide some sufficient conditions that guarantee local exponential stability under the assumptions posed in \Section~\ref{sec:problem_formulation}. Informally, we assume that each node updates its local variables and communicates with its neighbors infinitely often, and that the number of consecutive packet losses is bounded. Formally, we make the following assumptions. 
\begin{assumption}[Communications are persistent]
	For any iteration $k \in \N $ there exists a positive integer number $\tau$ such that each node performs at least one broadcast transmission within the interval $[k, k+\tau]$, i.e., for each $i \in \until{N}$ there exists $h \in [k, k+\tau]$ such that node $i$ is activated at iteration $h$.	%
	\label{ass:persistent_communications}
\end{assumption}
\begin{assumption}[Packet losses are bounded]
	There exists a positive integer $L$ such that the number of consecutive communication failures over every directed edge in the communication graph is smaller than $L$.
	\label{ass:bounded_packet_losses}
\end{assumption}

From the above two assumptions, it follows that, given $i \in V$ and $j \in \mathcal{N}_i^{\text{out}}$, node $j$ receives information from node $i$ at least once within the interval $\left[k, k+L\tau \right]$.

We now want to characterize the convergence properties of the Asymmetric
broadcast \ac{raNRC} algorithm. 
To do so, we start by introducing two Lemmas which will be later used.

Let $\bm{x}=\left[x_1,\ldots, x_N\right]^T$ and
$\bm{x}^0=\left[x_1^0,\ldots, x_N^0\right]^T$.  In the first lemma we show that
if the variable $\bm{x}$ is kept constant, then the components of the vector
$\bm{p}$ achieve consensus to the ratio
$\overline h(x_1,\ldots,x_N) / \overline g(x_1,\ldots,x_N)$. Viceversa, in the
second lemma, we show that if the components of $\bm{p}$ have reached consensus,
then the vector $\bm{x}$ exponentially converges to the global minimizer.

Formally, to state the first result, for a given $\bar{k}$, we consider the following dynamics, for $k \geq \bar{k}$,
\begin{equation}\label{eq:tilde_xi_NR}
\bm{\xi}_{\bar k}(k+1)=\varphi \left(k,  \bm{x}(\bar{k}), \bm{\xi}_{\bar k}(k) \right),
\end{equation}
initialized by $\bm{\xi}_{\bar k}(\bar k)= \bm{\xi}(\bar k)$. Observe that, $\bm{\xi}_{\bar k}$ describes the evolution of the variable $\bm{\xi}$, starting at iteration $\bar{k}$, assuming that the variable $\bm{x}$ is kept constant for $k \geq \bar{k}$, that is, $\bm{x}(k)= \bm{x}(\bar{k})$ for all $k \geq \bar{k}$.
In particular, in this scenario, we are interested in the behavior of the variable $\bm{p}$, that is, of the last block of components of $\bm{\xi}_{\bar{k}}$, that, in this case, similarly to $\bm{\xi}_{\bar k}$, we denote as $\bm{p}_{\bar{k}}$. We have the following result.

\begin{lemma}\label{lem:FirstResult}
For a given $\bar{k}$, consider, for $k \geq \bar{k}$, the dynamics in \eqref{eq:tilde_xi_NR}. Then, under \Assumptions~\ref{ass:persistent_communications},~\ref{ass:bounded_packet_losses}, we have that the point 
{
$$
\frac{\sum_\ell \bm{g}_\ell(\bm{x}_\ell(\bar{k}))}{\sum_\ell \bm{h}_\ell(\bm{x}_\ell(\bar{k}))}  \1
$$}
is exponentially stable for the variable $\bm{p}_{\bar{k}}$, that is, defined
{
$$
\tilde{\bm{p}}_{\bar{k}}(k):=  \bm{p}_{\bar{k}}(k) - \frac{\sum_\ell \bm{g}_\ell(\bm{x}_\ell(\bar{k}))}{\sum_\ell \bm{h}_\ell(\bm{x}_\ell(\bar{k}))}  \1,
$$}
there exists $C_{\bar{k}}>0$ and $0 \leq \rho_{\bar{k}} < 1$ such that 
\begin{equation}
\| \tilde{\bm{p}}_{\bar k}(\,k\,) \| \leq \,C_{\bar{k}}\, \rho_{\bar{k}}^{k-\bar{k}}\,\, \| \tilde{\bm{p}}_{\bar k}(\,\bar{k}\,)    \|.
\end{equation}
\end{lemma}
\begin{IEEEproof}
In the following we denote by $\bm{y}_{a;\bar{k}}$, $\bm{z}_{a; \bar{k}}$ the block components of $\bm{\xi}_{\bar k}$ corresponding to $\bm{y}_{a}$, $\bm{z}_{a}$.
To study the evolution of $\bm{p}_{\bar{k}}(k)$, we analyze the behavior of the variables $\bm{y}_{a;\bar{k}}(k)$, $\bm{z}_{a;\bar{k}}(k)$, separately. Consider $\bm{y}_{a;\bar{k}}(k)$. Observe that, since $\bm{x}(k)=\bm{x}(\bar{k})$, $k \geq \bar{k}$, and according to Assumptions \ref{ass:persistent_communications} and \ref{ass:bounded_packet_losses}, we have that there exists $\bar{k}' > \bar{k}$ such that $\bm{g}(k)= \bm{g}^{\text{old}}(k)$ for all $k \geq \bar{k}'$ and, hence, 
$$
\bm{y}_{a, \bar{k}}(k+1)= M(k) \bm{y}_{a, \bar{k}}(k),
$$
for $k \geq \bar{k}'$. A similar reasoning holds for $\bm{z}_{a, \bar{k}}(k)$. It follows that the variables $\bm{y}_{a, \bar{k}}(k)$, $\bm{z}_{a, \bar{k}}(k)$ and, in turn, the variables $\bm{y}_{\bar{k}}(k)$, $\bm{z}_{\bar{k}}(k)$ 
run the iterations of a ratio-consensus algorithm for $k \geq \bar{k}'$, as described in \cite{IFAC17}. 

From Lemma \ref{eq:ConservationProperty}, we have that, for $k \geq \bar{k}$ and, in particular, for $k \geq \bar{k}'$
\begin{align*}
\1^T \bm{y}_{a, \bar{k}}(k) & =  \sum_{\ell=1}^N  g_\ell\left(x_\ell(\bar{k})\right)\\
\1^T \bm{z}_{a, \bar{k}} (k) & =  \sum_{\ell=1}^N  h_\ell\left(x_\ell(\bar{k})\right).
\end{align*}
From Theorem $3$ in \cite{IFAC17}, it follows that $\frac{\bm{y}_{\bar{k}}(k)}{{\bm{z}_{\bar{k}}(k)}}$ converges exponentially to {$\frac{\sum_\ell \bm{g}_\ell(\bm{x}_\ell(\bar{k}))}{\sum_\ell \bm{h}_\ell(\bm{x}_\ell(\bar{k}))}  \1$}. \\
\end{IEEEproof}
Now, let us assume that, for each $k$, the variable $\bm{p}(k)$ has reached consensus and consider the following dynamics for the variable $\bm{x}$, 
\begin{align}\label{eq:reduced_system}
\bm{x}(k+1)
&= \bm{x} (k)+ \varepsilon \, S(k) \, \left(\frac{\sum_\ell g_\ell(x_\ell(k))}{{\sum_\ell h_\ell(x_\ell(k))}}  \1 -\bm{x}(k) \right) \nonumber\\
&= \bm{x} (k) + \varepsilon\, \tilde{\phi} (k; \bm{x}(k))
\end{align}
where
$$
\tilde{\phi} (k; \bm{x}(k)) =  S(k) \, \left(\frac{\sum_\ell g_\ell(x_\ell(k))}{{\sum_\ell h_\ell(x_\ell(k))}}  \1 \nonumber-\bm{x}(k) \right).
$$
Let 
\begin{equation}\label{eq:N-optimum}
\bm{x}^* = x^* \1,
\end{equation}
where we recall that $x^*$ is the minimizer of the optimization problem in \eqref{eq:problem}. By standard algebraic manipulations, one can see that $\bm{x}^*$ is an equilibrium of  \eqref{eq:reduced_system}. The following result states that the linearized version of  \eqref{eq:reduced_system} around $\bm{x}^*$ is an exponentially stable system.

\begin{lemma}\label{lem:SecondResult}
Consider system in \eqref{eq:reduced_system} and let $\bm{x}^*$ be as in \eqref{eq:N-optimum}. Let 
$$
A(k)= I+ \epsilon \frac{\partial \tilde{\phi}}{\partial \bm{x}}(k; \bm{x}) |_{\bm{x}=\bm{x}^*},
$$ 
and, accordingly, consider the auxiliary system
\begin{equation}\label{eq:Linearized_NR}
\tilde{\bm{x}}(k+1)=A(k) \tilde{\bm{x}}(k).
\end{equation}
Then, under \Assumptions~\ref{ass:persistent_communications},~\ref{ass:bounded_packet_losses}, $\tilde{\bm{x}}=0$ is exponentially stable equilibrium point for \eqref{eq:Linearized_NR}.
\end{lemma}
\begin{IEEEproof}
Let 
$$
\alpha(\bm{x}(k))= \frac{\sum_\ell g_\ell(x_\ell(k))}{\sum_\ell h_\ell(x_\ell(k))}.
$$
Computing the partial derivative of $\alpha$ with the respect to $x_i$ we get
\begin{align*}
&\left[\frac{\partial \alpha}{\partial x_i} \right] |_{\bm{x}=\bm{x}^*}= \frac{g'_i(x^*) \,\sum_{\ell=1}^N  h_\ell \left(x^*\right) - h'_i\left(x^*\right) \, \sum_{\ell=1}^N  g_\ell\left(x^*\right)}{\left(\sum_{\ell=1}^N  h_\ell \left(x^*\right)\right)^2}
\end{align*}
with
\begin{align*}
&g'_i(x^*) \,\sum_{\ell=1}^N  h_\ell \left(x^*\right)\, -\, h'_i\left(x^*\right) \, \sum_{\ell=1}^N  g_\ell\left(x^*\right)\\
&= \left(\,f'''_i (x^*) \,x^*\, +f''_i(x^*) \,-\, f''_i(x^*) \,\right) \,\sum_{\ell=1}^N  f''_\ell \left(x^*\right) \\
&\qquad -f'''_i(x^*) \, \sum_{\ell=1}^N \left( f''_\ell\left(x^*\right)\,x^*\, -f'_\ell(x^*)\right)\\
& = f'''_i (x^*)\, x^* \, \sum_{\ell=1}^N  f''_\ell \left(x^*\right) \, - \, f'''_i(x^*) \,x^*\, \sum_{\ell=1}^N  f''_\ell\left(x^*\right)\\
& \qquad + f'''_i(x^*) \, \sum_{\ell=1}^N f'_\ell(x^*)\\
&=0,
\end{align*}
where , in the last equality, we have used the fact that $\sum_{\ell=1}^N f'_\ell(x^*)=0$. From the previous calculations, it turns out that 
$$
A(k)=I -\varepsilon S(k).
$$
By Assumption \ref{ass:persistent_communications}, we have that the matrix
$$
\bar{A}_{k,\tau}=\prod_{s=k}^{k+\tau} A(k),
$$
is a diagonal matrix such that $0<\left[\bar{A}_{k,\tau}\right]_{ii} <1-\varepsilon$, for all $i$. Then, system in \eqref{eq:reduced_system} satisfies the stated property.
\end{IEEEproof}


Intuitively, one would conclude that when the parameter $\varepsilon$ is small
the results of the two lemma can be combined to simultaneously obtain asymptotic
consensus and convergence to the global minimizer. This is formally shown in the
next theorem which characterizes the convergence properties of the Asymmetric
broadcast \ac{raNRC} algorithm.%

\begin{theorem}
	Under \Assumptions~\ref{ass:persistent_communications},~\ref{ass:bounded_packet_losses} and the assumptions posed in \Section~\ref{sec:problem_formulation}, there exist some positive scalars $\varepsilon_c$ and $\delta$ such that, if the initial conditions $\bm{x}^0 \in \R^N$ satisfy $\|\bm{x}^o-x^*\1\|<\delta$ and if $\varepsilon$ satisfies $0 < \varepsilon < \varepsilon_c$ then the local variables $x_i$ in \Algorithm~\ref{alg:raNRC} are exponentially stable with respect to the global minimizer $x^*$.
	\label{thm:convergence}
\end{theorem}

\begin{IEEEproof}
The proof of the result is based on showing that the system in \eqref{eq:NonLinear_System} satisfies the assumptions of Proposition \ref{prop:MainPropApp}. To do so, we start by defining, for $k \geq \bar{k}$,
\begin{equation}\label{eq:xi*_NR}
\bm{\xi}^*_{ \bm{x}(\bar{k}), \bm{\xi}(\bar{k})} (k)=  \tilde{I} \, \bm{\xi}_{\bar k}(k) \, + \, \tilde{u},
\end{equation}
where $\bm{\xi}_{\bar k}(k)$ is defined as in \eqref{eq:tilde_xi_NR} and 
where
$$
 \tilde{I} = 
 \left[ 
\begin{array}{cc}
I_{(6N+2|\mathcal{E}|) \times (6N+2|\mathcal{E}|)} & 0_{(6N+2|\mathcal{E}|) \times N} \\
0_{N \times (6N+2|\mathcal{E}|)} & 0_{N \times N}
\end{array}
\right]
$$
and
$$
\tilde{u}=
\left[ 
\begin{array}{c}
0_{(6N+2|\mathcal{E}|) \times 1} \\
\bm{p}_{\bm{x}(\bar{k})}^*
\end{array}
\right],
$$
with 
$$
\bm{p}^*_{\bm{x}(\bar{k})}=   \frac{\sum_\ell \bm{g}_\ell(x_\ell(\bar{k}))}{\sum_\ell \bm{h}_\ell(x_\ell(\bar{k}))}  \1.
$$
Observe that the first six blocks components of $\bm{\xi}^*_{ \bm{x}(\bar{k}), \bm{\xi}(\bar{k})} (k)$ coincide with the first six blocks components of $\bm{\xi}_{\bar k}(k)$, while the last block component is constant for all $k\geq \bar{k}$.
Moreover, for $k \geq \bar{k}$, let 
$$
\tilde{\bm{\xi}}_{\bar k}(k):=\bm{\xi}_{\bar k}\left(k\right)-\bm{\xi}_{\bm{x}(\bar{k}), \bm{\xi}(\bar{k})}^*\left(k\right).
$$
Based on the previous observation, we have that the first six blocks components of $\tilde{\bm{\xi}}_{\bar k}(k)$ are equal to zero, while the last block component is equal to 
$$
\bm{p}_{\bar{k}}(k)- \bm{p}^*_{\bm{x}(\bar{k})}.
$$
From Lemma \ref{lem:FirstResult}, it easily follows that there exists $C_{\bar{k}}>0$ and $0 \leq \rho_{\bar{k}} < 1$ such that
\begin{equation}
\|\tilde{\bm{\xi}}_{\bar k}(\,k\,)\| \leq \,C_{\bar{k}}\, \rho_{\bar{k}}^{k-\bar{k}}\,\, \| \tilde{\bm{\xi}}_{\bar k}(\,\bar{k}\,)\|.
\end{equation}
This shows that system in \eqref{eq:NonLinear_System} satisfies property in \eqref{eq:ex_xi}, in Appendix \ref{app:nonlinear}

Consider now the system
\begin{align}\label{eq:reduced_system_1}
\bm{x}(k+1)&= \bm{x}(k)+\varepsilon \, \phi \left(k, \bm{x}(k),\bm{\xi}_{\bm{x}(k),\bm{\xi}(k)}^*(k) \right)\\
&= \bm{x} (k)+ \varepsilon \, S(k) \, \left(\frac{\sum_\ell g_\ell(x_\ell(k))}{\sum_\ell h_\ell(x_\ell(k))}  \1 -\bm{x}(k) \right) \nonumber\\
&= \bm{x} (k) + \varepsilon\, \tilde{\phi} (k; \bm{x}(k))
\end{align}

In Lemma \ref{lem:SecondResult}, it is established that the previous system satisfies Assumption \ref{ass:Linearized}, in Appendix \ref{app:nonlinear}.
Hence, Proposition \ref{prop:MainPropApp}, in Appendix \ref{app:nonlinear} can be applied to system in \eqref{eq:NonLinear_System}, yielding the result of the statement.
\end{IEEEproof}

\begin{remark}
	\Algorithm~\ref{alg:raNRC} assumes the initial conditions of the local variable $x_i$ to be all identical to $x^o$. Although not being a very stringent requirement, this assumption can be relaxed, that is, slightly modified versions of \Theorem~\ref{thm:convergence} would hold even in the case $x_{i} = x_{i}^{o}$ as soon as all the initial conditions are sufficiently close to the global minimizer $x^*$, i.e., as soon as $|x_i^o-x^*|<\delta$ for all $i = 1, \ldots, N$.
\end{remark}

\begin{remark}
The initial conditions on the local variables $y_i = g_i^{\text{old}}=g_i=f''_i(x^o)x^o-f'_i(x^o)$ and $z_i = h_i^{\text{old}}=h_i=f''_i(x^o)$ are instead more critical for the convergence of the local variables $x_i$ to the true minimizer $x^*$. As shown in~\cite{zanella2011newton}, small perturbations of these initial conditions can lead to convergence to a point $\overline{x} \neq x^*$ (notice that these perturbations do not affect the stability of the algorithm, so that possible small numerical errors due to the computation and data quantization do not disrupt the convergence properties of the algorithm). Moreover, the map from the amplitude of these perturbations and the distance $\left\| \overline{x} - x^* \right\|$ is continuous, so that if these perturbations are small then $\overline{x} \approx x^*$.
\end{remark}

\begin{remark}
Although the previous theorem guarantees only local exponential convergence, numerical simulations on real datasets seem to indicate that the basin of attraction is rather large and stability is mostly dictated by the choice of the parameter $\varepsilon$. { However, for the special but relevant case when the cost functions $f_i(x)$ are quadratic, as in distributed least-squares problems, local stability implies global stability~\cite{ECC15}.}
\end{remark}

\begin{remark}
The major challenges in proving the main results are related to proving that the ra-NRC algorithm satisfy a number of technical conditions required by standard theory of separation of time-scales. Different conditions and theorems are available for continuous time dynamical systems (we refer the interested reader to Chapter~11 in \cite{khalil1996noninear}). In particular, we are interested in proving exponential stability for a non-autonomous discrete time dynamical system whose closest counterpart in the continuous time is given by Theorem~11.4 in \cite{khalil1996noninear}. Besides some standard conditions on smoothness and uniformity of the dynamical flows involved, there are three major requirements that need to be satisfied: the first is that the fast dynamics converges exponentially to an equilibrium manifold, the second is that the slow dynamics restricted to this manifold is exponentially stable, and the third is that a number of \emph{bounded interconnection conditions} which represent the perturbation of the slow dynamics into the fast dynamics and vice-versa, are satisfied. As for the first requirement, we were able to guaranteed it by extending (see \cite{IFAC17}) the work by \cite{vaidya2011}, which only provided convergence in probability. As for the second one, we are able to prove local exponential stability of the slow dynamics which is not trivial since the dynamics is non-autonomous. As for the last requirement on the \emph{bounded interconnection conditions}, very much depends on cost functions and in the discrete-time domain it is difficult to provide \emph{global} guarantees. However, under some mild smoothness conditions, we were able to show that the conditions on \emph{bounded interconnection conditions} are locally satisfied, and, in turn, to prove local exponential stability. 
\end{remark}

\section{Numerical Experiments}
\label{sec:numerical_experiments}

\begin{figure}[!htbp]
	\begin{center}
	\input{figure_network_r05}
	\end{center}
	\caption{The random geometric network considered in the simulations.}
	\label{fig:networks}
\end{figure}
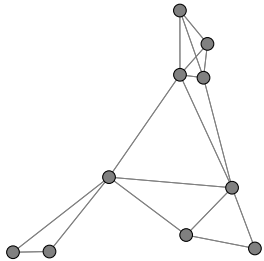

We consider a random geometric network with $10$ nodes in $[0, 1]^{2}$ and with communication radius $r = 0.5$ as in \Figure~\ref{fig:networks}.

As cost functions, we consider the distributed training of a Binomial-Deviance
based spam-nonspam
classifier{~\cite[Chap.~10.5]{friedman2001elements} } 
where the training set is a database of $E$ emails with $j$ the email index, $y_{j} = -1, 1$ indicating if email $j$ is spam or not, $\chi_{j} \in \mathbb{R}^{n-1}$ summarizing the $n-1$ features of the $j$-th email (in our case the frequency of words ``make'', ``address'', and ``all''). Letting $x=\left( x', x_{0} \right) \in \mathbb{R}^{n-1} \times \mathbb{R}$ represent a generic classification hyperplane, distributedly training a Binomial-Deviance based classifier corresponds to solve the distributed optimization problem with local costs defined by
\begin{equation}
    f_i\left( x \right)
    \DefinedAs
    \sum_{j \in E_i}
    \log
    \Big(
        1 
        +
        \exp
        \left(
			- y_{j} \left( \chi_{j}^T x' + x_{0} \right)
        \right)
    \Big)
    +
    \gamma \left\| x' \right\|^{2}_{2}
    \label{equ:spam_nonspam_classification}
\end{equation}
where $E_i$ is the set of emails available to agent $i$, $E=\cup_{i=1}^N E_i$, and $\gamma$ is a global regularization parameter. In our experiments we consider $|E| = 5000$ emails from the spam-nonspam UCI repository\footnote{\texttt{http://archive.ics.uci.edu/ml/datasets/Spambase}}, randomly assigned to the 10 nodes users communicating as in \Figure~\ref{fig:networks}.  As a performance index, we consider the Mean Squared Error (MSE)
$$
    \frac{1}{N}
    \sum_{i = 1}^{N} 
    \left\| x_{i}(k) - x^{*} \right\|^{2}
$$
as a function of the iteration index $x$.

\begin{figure}[!htbp]
	\centering
	\input{figure_time_evolution_varying_varepsilon}
	\caption{Evolution of the MSE in time for a typical realization of the optimization process for a fixed packet loss probability of $0.1$ and different values of $\varepsilon$.}
	\label{fig:time-evolution}
\end{figure}
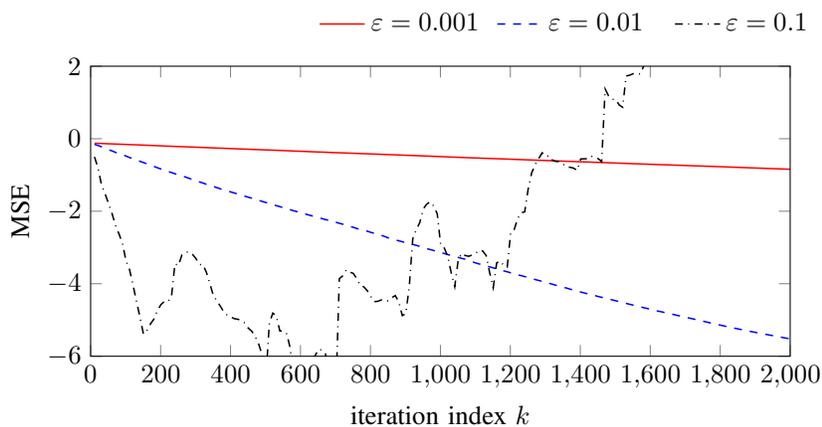

\Figure~\ref{fig:time-evolution} then plots the evolution of the MSE of a typical realization of the optimization process as a function of the iteration index $k$, for a fixed packet loss probability equal to $0.1$, and for different values of $\varepsilon$. The figure confirms the intuition that increasing $\varepsilon$ may lead to faster convergence properties, but only up to a certain value; too aggressive $\varepsilon$'s may indeed hinder the convergence property of the algorithm.

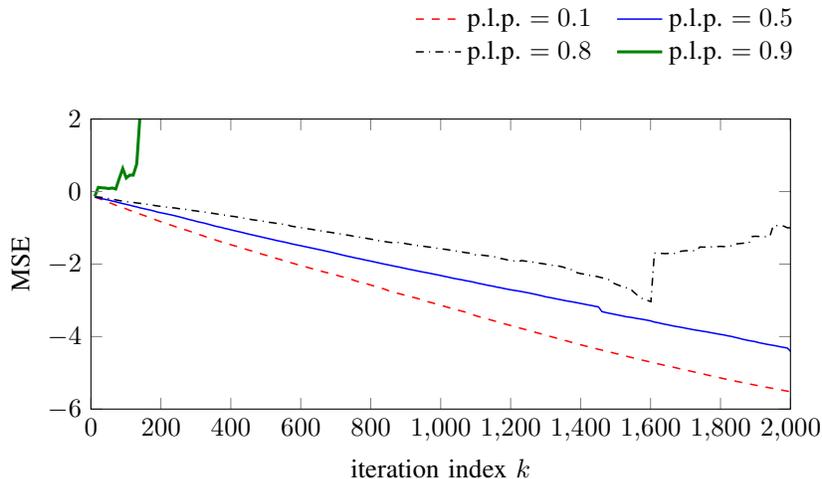
\begin{figure}[!htbp]
	\centering
	\input{figure_time_evolution_varying_plp}
	\caption{Effect of varying the probability of packets losses on the MSE of single realizations of the optimization process for $\varepsilon = 0.01$.}
	\label{fig:fixed-varepsilon}
\end{figure}

\Figure~\ref{fig:fixed-varepsilon} instead inspects the effect of varying the probability of packets losses on the MSEs of single realizations for a fixed $\varepsilon$. This figure confirms the intuition that, independently of $\varepsilon$, increasing the chances of packet losses leads to initially slower convergence properties and eventually divergent behaviors.

\section{Conclusions}
\label{sec:conclusions}

Implementations of distributed optimization methods in real-world scenarios require strategies that are both able to cope with real-world problematics (like unreliable, asynchronous and directed communications), and converge sufficiently fast so to produce usable results in meaningful times. Here we worked towards this direction, and improved an already existing distributed optimization strategy, previously shown to have fast convergence properties, so to make it tolerate the previously mentioned real-world problematics.

More specifically, we considered a robustified version of the Newton-Raphson consensus algorithm originally proposed in~\cite{zanella2011newton} and proved its convergence properties under some general mild assumptions on the local costs. From technical perspectives we shown that under suitable assumptions on the initial conditions, on the step-size parameter, on the connectivity of the communication graph and on the boundedness of the number of consecutive packet losses, the considered optimization strategy is locally exponentially stable around the global optimum as soon as the local costs are $\mathcal{C}^{2}$ and strongly convex with second derivative bounded from below.

We also shown how the strategy can be applied to real world scenarios and datasets, and be used to successfully compute optima in a distributed way.

We then notice that the results offered in this manuscript do not deplete the set of open questions and plausible extensions of the Newton Raphson consensus strategy. We indeed devise that the algorithm is potentially usable as a building block for distributed interior point methods, but that some lacking features prevent this development. Indeed it is still not clear how to tune the parameter $\varepsilon$ online so that the convergence speed is dynamically adjusted (and maximized), how to account for equality constraints of the form $Ax = b$, and how to update the local variables $x_{i}$ using partition-based approaches so that each agent keeps and updates only a subset of the components of $x$.

\appendices
\section{General results on discrete-time nonlinear systems}\label{app:nonlinear} 
The proofs and results of this appendix can be found in the technical report \cite{LyapunovDiscreto}.

\noindent
Consider the system 
\begin{equation}\label{eq:MainSystem}
\left\{
\begin{array}{l}
\bm{x}(k+1)=\bm{x}(k)+\varepsilon \phi (k, \bm{x}(k), \bm{\xi}(k))\\
\bm{\xi}(k+1)=\varphi (k, \bm{x}(k)), \bm{\xi}(k))
\end{array}
\right.
\end{equation}
where $\bm{x} \in \R^{n_1}$, $\bm{\xi} \in \R^{n_2}$, $\phi \,:\, \N \times \R^{n_1} \times \R^{n_2}\,\to\, \R^{n_1}$, $\varphi \,:\, \N \times \R^{n_1} \times \R^{n_2}\,\to\, \R^{n_2}$, $\varepsilon>0$ and with given initial conditions $\bm{x}(0)$, $\bm{\xi}(0)$.

\noindent
For a given $\bar{k} \in\N$, consider the system, for $k \geq \bar{k}$, 
\begin{equation}\label{eq:tilde_xi}
\bm{\xi}_{\bar k}(k+1)=\varphi \left(k, \,\bm{x}(\bar{k}), \,\tilde{\bm{\xi}}_{\bar k}(k) \right),
\end{equation}
initialized by $\bm{\xi}_{\bar k}(\bar k)= \bm{\xi}(\bar k)$, where $\bm{\xi}(\bar k)$ is obtained ruling system \eqref{eq:MainSystem} up to $\bar{k}$.\\
Given $\bar{k}$, assume that, for $k \geq \bar{k}$, there exists a sequence
\begin{equation}\label{eq:sequence}
k \, \to \,\bm{\xi}_{\bm{x}(\bar{k}), \bm{\xi}(\bar{k})}^*\left(k\right),
\end{equation}
in general dependent on $\bm{x}(\bar{k})$ and $\bm{\xi}(\bar{k})$,
such that the evolution
\begin{equation}\label{eq:xi*}
\tilde{\bm{\xi}}_{\bar k}(k):=\bm{\xi}_{\bar{k}}\left(k\right)-\bm{\xi}_{\bm{x}(\bar{k}), \bm{\xi}(\bar{k})}^*\left(k\right)
\end{equation}
satisfies the property
\begin{equation}\label{eq:ex_xi}
\|\tilde{\bm{\xi}}_{\bar k}(k)\| \leq \, C_{\bar{k}}\, \rho_{\bar{k}}^{k-\bar{k}}\,\, \| \tilde{\bm{\xi}}_{\bar k}(\bar{k})\|,
\end{equation}
for suitable $C_{\bar{k}}>0$ and $0 \leq \rho_{\bar{k}} < 1$, that is $\tilde{\bm{\xi}}'_{\bar k}=0$ is an exponentially stable point for the evolution in \eqref{eq:xi*}. 
Basically, the property in \eqref{eq:ex_xi} establishes that there exists a trajectory $\bm{\xi}^*$ to which the trajectory of the variable $\bm{\xi}$, generated keeping the variable $\bm{x}$ constant, converges asymptotically.  

Next, let us assume that, for each $k$, the variable $\bm{\xi}$ has already reached the asymptotic convergence to the corresponding trajectory $\bm{\xi}^*$. More precisely, observe that there exists a family of sequences of the type \eqref{eq:sequence}, where each sequence starts from a different index $k$. From this family we can build the following new sequence 
\begin{equation}\label{eq:new_sequence}
k \, \to \,\bm{\xi}_{\bm{x}(k), \bm{\xi}(k)}^*\left(k\right),
\end{equation}
where, to the index $k$, we have associated the first element of the sequence which starts at $k$. Based on \eqref{eq:new_sequence},
we consider the system 
\begin{equation}\label{eq:system1}
\bm{x}(k+1)=\bm{x}(k)+\varepsilon \, \phi \left(k, \bm{x}(k),\bm{\xi}_{\bm{x}(k),\bm{\xi}(k)}^*(k) \right).
\end{equation}
Assume that $\bm{\xi}_{\bm{x}(k),\bm{\xi}(k)}^*(k)$ is such that there exists a suitable map $\tilde{\phi}: \N \times \R^{n_1} \,\to\, \R^{n_1}$ such that
\eqref{eq:system1} can be, equivalently, rewritten as
\begin{equation}\label{eq:system1-tilde}
\bm{x}(k+1)=\bm{x}(k)+\varepsilon  \, \tilde{\phi} \left(k, \bm{x}(k) \right),
\end{equation}
that is, $\tilde{\phi} \left(k, \bm{x}(k) \right)=  \phi \left(k, \bm{x}(k),\bm{\xi}_{\bm{x}(k),\bm{\xi}(k)}^*(k) \right)$.
We make the following assumption.
\begin{assumption}\label{ass:Linearized}
Let $\bm{x}^*$ be an equilibrium point for \eqref{eq:system1-tilde}. We assume that, there exists $r>0$ such that $\tilde{\phi}$ is continuously differentiable on $D=\left\{\bm{x} \in \R^{n_1} \,| \, \|\bm{x}-\bm{x}^*\| < r \right\}$ and the Jacobian matrix $[\partial \tilde{\phi} / \partial \bm{x}]$ is bounded and Lipschitz on $D$, uniformly in $k$. In addition, defining 
$$
A(k)= I+ \varepsilon \frac{\partial \tilde{\phi}}{\partial \bm{x}}(k; \bm{x}) |_{\bm{x}=\bm{x}^*},
$$ 
and considering the auxiliary system
\begin{equation}\label{eq:A(t)}
\tilde{\bm{x}}(k+1)=A(k) \tilde{\bm{x}}(k),
\end{equation}
we assume that $\tilde{\bm{x}}=0$ is exponentially stable equilibrium point for \eqref{eq:A(t)}.
\end{assumption}

The following Proposition characterizes the convergence properties of system \eqref{eq:MainSystem}.

\begin{proposition}\label{prop:MainPropApp}
Consider system in \eqref{eq:MainSystem}. For any $\bar{k}$, assume that there exists a sequence as in \eqref{eq:sequence} such that property \eqref{eq:ex_xi} is satisfied. Consider system in \eqref{eq:system1}. Let $\bm{x}^*$ be an equilibrium point for \eqref{eq:system1}. Assume Assumption \eqref{ass:Linearized} holds true. Then, there exist $r>0$ and $\varepsilon^*>0$, such that, for all $\varepsilon \in (0, \varepsilon^*]$ and for all $x(0) \in B_r^n=\left\{\bm{x} \in \R^n : \|\bm{x}-\bm{x}^*\| < r\right\}$, the trajectory $\bm{x}(t)$ generated by \eqref{eq:MainSystem}, 
converges exponentially to $\bm{x}^*$, i.e., there exist $C>0$ and $0<\lambda<1$ such that
$$
\|\bm{x}(k)-\bm{x}^*\| \leq C \lambda ^k \|\bm{x}(0)-\bm{x}^*\|.
$$
\end{proposition}

\addcontentsline{toc}{section}{References}
\bibliographystyle{IEEEtran}
\bibliography{bibliography} 

\end{document}

%% file: preamble.tex
\usepackage{amsmath,amssymb,amsfonts,theorem}
\usepackage{graphicx,epsfig,color}
\usepackage{subfigure,psfrag}
\usepackage{url}
\usepackage{bm}
\usepackage{acronym}
\usepackage[color=white!90!black, textwidth=3cm]{todonotes}
\usepackage[switch, pagewise, mathlines]{lineno}
\usepackage{algorithm}												%
\usepackage{algpseudocode}											
\usepackage{pgfplots}
\usepackage[shortlabels]{enumitem} 
\usepackage{algcompatible}
\usepackage{soul}

\graphicspath{{../Images/}}

\newcommand{\until}[1]{\{1,\dots, #1\}}

\def\1{\mathbf{1}}

\newcommand\oprocendsymbol{\hbox{$\square$}}
\newcommand\oprocend{\relax\ifmmode\else\unskip\hfill\fi\oprocendsymbol}

\def\V{\mathcal{V}}
\def\E{\mathcal{E}}
\def\C{\mathcal{C}}

\def\N{\mathbb{N}}
\def\R{\mathbb{R}}

\newcommand{\argmin}{\mathop{\operatorname{argmin}}}

\newtheorem{theorem}{Theorem}[section]
\newtheorem{proposition}[theorem]{Proposition}
\newtheorem{corollary}[theorem]{Corollary}

\newtheorem{lemma}[theorem]{Lemma}
\newtheorem{assumption}[theorem]{Assumption}

{\theorembodyfont{\rmfamily} 
\newtheorem{remark}[theorem]{Remark}

}

\def\log{\mathop{\operatorname{log}}}

\newcommand	{\Section}				[0]	{Section}

\newcommand	{\Figure}				[0]	{Figure}

\newcommand	{\Theorem}				[0]	{Theorem}

\newcommand	{\Algorithm}			[0]	{Algorithm}

\newcommand	{\Assumptions}			[0]	{Assumptions}

\newcommand{\DefinedAs}				[0]	{\mathrel{\mathop:}=}
\newcommand{\IDefinedAs}			[0]	{=\mathrel{\mathop:}}

\newcommand{\cmax}					[1]	{\left[ #1 \right]_{c}}

\ifarXiv
	\newcommand{\GN}[1]{{#1}}
	
\else
	\newcommand{\GN}[1]{{\color{blue} #1}}
	
\fi

%% file: figure_NR_lossy_D.tex
\tikzset
{
	NRNode/.style =
	{
		draw,
		rectangle,
		minimum width	= 2.0cm,
		minimum height	= 0.5cm,
		rounded corners	= 0.2cm,
		fill			= black!10!white,
		align			= center,
		thick,
	},
	Dots/.style =
	{
		line width		= 3pt,
		line cap		= round,
		dash pattern	= on 0pt off 2\pgflinewidth,
		shorten <		= 0.12cm,
		shorten >		= 0.0cm
	},
	Arrow/.style =
	{
		-latex,
		rounded corners	= 0.2cm,
		thick,
	},
	IArrow/.style =
	{
		Arrow,
		latex-,
	},
}

\begin{tikzpicture}

	\node (NR1) [NRNode] {\acl{NR} \\ $x_{1}$};
	\node (NRN) [NRNode, below = 0.9cm of NR1] {\acl{NR} \\ $x_{N}$};
	\draw [Dots] (NR1) -- (NRN);
	\node (Coo) [NRNode, below = 0.9cm of NRN] {Robust Asynchronous \\ Ratio Consensus \vspace{0.2cm} \\ $\sigma_{i,y}, \rho_{i,y}^{(j)}, \sigma_{i,z}, \rho_{i,z}^{(j)}$};

	\node [above = 0.1cm of NR1, text = red!80!black] {\emph{local computation}};
	\node [below = 0.1cm of Coo, text = red!80!black] {\emph{local cooperation}};

	\draw [Arrow] (NR1) -- ++(3.2cm,0) node [pos = 0.0, above right] {$g_{1}, h_{1}$} |- (Coo.345);
	\draw [Arrow] (NRN) -- ++(2.7cm,0) node [pos = 0.0, above right] {$g_{N}, h_{N}$} |- (Coo.15);

	\draw [IArrow] (NR1) -- ++(-3.2cm,0) node [pos = 0.0, above left] {$y_{1}, z_{1}$} |- (Coo.195);
	\draw [IArrow] (NRN) -- ++(-2.7cm,0) node [pos = 0.0, above left] {$y_{N}, z_{N}$} |- (Coo.165);

\end{tikzpicture}

%% file: figure_network_r05.tex
\begin{tikzpicture}
[
	xscale	= 1,	
	yscale	= 1,	
]
	\ifarXiv
		\begin{axis}
		[
			width 	= 0.3\columnwidth,
			height	= 0.3\columnwidth,
			hide x axis,
			hide y axis,
		]
	\else
		\begin{axis}
		[
			width 	= 0.6\columnwidth,
			height	= 0.6\columnwidth,
			hide x axis,
			hide y axis,
		]
	\fi
		\addplot
		[
			only marks,
			mark options	=
			{
				scale 		= 1.2,
				draw		= black,
				fill		= black!50!white,
				opacity		= 1
			},
		]
		table
		{network.txt.NodesPositions};
		\addplot
		[
			quiver		=
			{
				u		= \thisrow{towardsx},
				v		= \thisrow{towardsy}
			},
			solid,
			color		= black!50!white,
			line width	= 0.01cm,
		]
		table
		{network.txt.LinksList};
	\end{axis}
\end{tikzpicture}

%% file: figure_time_evolution_varying_varepsilon.tex
\begin{tikzpicture}
\ifarXiv
	\begin{axis}
	[
		width					= 0.6\columnwidth,
		height					= 0.3\columnwidth,
		axis on top,
		xmin					= 0,
		xmax					= 2000,
		ymax					= 2,
		ymin					= -6,
		xlabel					= {iteration index $k$},
		xlabel near ticks,
		ylabel					= {MSE},
		ylabel near ticks,
		legend style			=
		{
			draw				= none,
			fill				= white,
			at					= {(1.10,1.22)},
			nodes				= {text width = 1.6cm},
			cells				= {anchor = east},
		},
		legend columns			= -1, 
	]
\else
	\begin{axis}
	[
		width					= 0.9\columnwidth,
		height					= 0.5\columnwidth,
		axis on top,
		xmin					= 0,
		xmax					= 2000,
		ymax					= 2,
		ymin					= -6,
		xlabel					= {iteration index $k$},
		xlabel near ticks,
		ylabel					= {MSE},
		ylabel near ticks,
		legend style			=
		{
			draw				= none,
			fill				= white,
			at					= {(1.10,1.22)},
			nodes				= {text width = 1.6cm},
			cells				= {anchor = east},
		},
		legend columns			= -1, 
	]
\fi
		\addplot
		[
			solid,
			red, 
			line width = 0.02cm,
		]
		table [x = iteration, y = log_error]
		{NRC__N10_r0.5_costbinomialDeviance_epsilon0.001_packetLoss0.1.txt};
		\addlegendentry{$\varepsilon = 0.001$};
		\addplot
		[
			dashed,
			blue,
			line width = 0.02cm,
		]
		table [x = iteration, y = log_error]
		{NRC__N10_r0.5_costbinomialDeviance_epsilon0.01_packetLoss0.1.txt};
		\addlegendentry{$\varepsilon = 0.01$};
		\addplot
		[
			dashdotted,
			black,
			line width = 0.02cm,
		]
		table [x = iteration, y = log_error]
		{NRC__N10_r0.5_costbinomialDeviance_epsilon0.1_packetLoss0.1.txt};
		\addlegendentry{$\varepsilon = 0.1$};
	\end{axis}
\end{tikzpicture}

%% file: figure_time_evolution_varying_plp.tex
\begin{tikzpicture}
\ifarXiv
	\begin{axis}
	[
		width					= 0.6\columnwidth,
		height					= 0.3\columnwidth,
		axis on top,
		xmin					= 0,
		xmax					= 2000,
		ymax					= 2,
		ymin					= -6,
		xlabel					= {iteration index $k$},
		xlabel near ticks,
		ylabel					= {MSE},
		ylabel near ticks,
		legend style			=
		{
			draw				= none,
			fill				= white,
			at					= {(1.05,1.42)},
			nodes				= {text width = 1.9cm},
			cells				= {anchor = east},
		},
		legend columns			= 2, 
	]
\else
	\begin{axis}
	[
		width					= 0.9\columnwidth,
		height					= 0.5\columnwidth,
		axis on top,
		xmin					= 0,
		xmax					= 2000,
		ymax					= 2,
		ymin					= -6,
		xlabel					= {iteration index $k$},
		xlabel near ticks,
		ylabel					= {MSE},
		ylabel near ticks,
		legend style			=
		{
			draw				= none,
			fill				= white,
			at					= {(1.05,1.42)},
			nodes				= {text width = 1.9cm},
			cells				= {anchor = east},
		},
		legend columns			= 2, 
	]
\fi
		\addplot
		[
			dashed,
			red,
			line width = 0.02cm,
		]
		table [x = iteration, y = log_error]
		{NRC__N10_r0.5_costbinomialDeviance_epsilon0.01_packetLoss0.1.txt};
		\addlegendentry{p.l.p.\ $= 0.1$};
		\addplot
		[
			solid,
			blue,
			line width = 0.02cm,
		]
		table [x = iteration, y = log_error]
		{NRC__N10_r0.5_costbinomialDeviance_epsilon0.01_packetLoss0.5.txt};
		\addlegendentry{p.l.p.\ $= 0.5$};
		\addplot
		[
			dashdotted,
			black,
			line width = 0.02cm,
		]
		table [x = iteration, y = log_error]
		{NRC__N10_r0.5_costbinomialDeviance_epsilon0.01_packetLoss0.8.txt};
		\addlegendentry{p.l.p.\ $= 0.8$};
		\addplot
		[
			solid,
			line width = 0.04cm,
			green!50!black,
		]
		table [x = iteration, y = log_error]
		{NRC__N10_r0.5_costbinomialDeviance_epsilon0.01_packetLoss0.9.txt};
		\addlegendentry{p.l.p.\ $= 0.9$};
	\end{axis}
\end{tikzpicture}